\documentclass[12pt]{amsart}
\usepackage{amsfonts}
\usepackage{amsmath}
\usepackage{amssymb,latexsym}
\usepackage[mathcal]{eucal}
\usepackage{amscd}
\input xy
\xyoption{all}

\DeclareMathOperator*{\esssup}{ess-sup}

\newcommand{\Acal}{\mathcal{A}}
\newcommand{\Bcal}{\mathcal{B}}

\newcommand{\Dcal}{\mathcal{D}}
\newcommand{\Ecal}{\mathcal{E}}
\newcommand{\Fcal}{\mathcal{F}}

\newcommand{\Ical}{\mathcal{I}}

\newcommand{\Ncal}{\mathcal{N}}

\newcommand{\Pcal}{\mathcal{P}}
\newcommand{\Qcal}{\mathcal{Q}}
\newcommand{\Rcal}{\mathcal{R}}
\newcommand{\Scal}{\mathcal{S}}

\newcommand{\Ucal}{\mathcal{U}}
\newcommand{\Vcal}{\mathcal{V}}

\newcommand{\Xcal}{\mathcal{X}}
\newcommand{\Ycal}{\mathcal{Y}}
\newcommand{\Zcal}{\mathcal{Z}}

\newcommand{\ch}{\mathbf{1}}

\newcommand{\Lf}{\mathfrak{L}}

\newcommand{\Z}{\mathbb{Z}}

\newcommand{\R}{\mathbb{R}}
\newcommand{\C}{\mathbb{C}}
\newcommand{\N}{\mathbb{N}}
\newcommand{\T}{\mathbb{T}}

\newcommand{\A}{\mathbb{A}}

\newcommand{\Eb}{\mathbf{E}}

\newcommand{\Pb}{\mathbf{P}}
\newcommand{\Qb}{\mathbf{Q}}

\newcommand{\Wb}{\mathbf{W}}
\newcommand{\Xb}{\mathbf{X}}
\newcommand{\Yb}{\mathbf{Y}}
\newcommand{\Zb}{\mathbf{Z}}

\newcommand{\XT}{$(X,T)$\ }
\newcommand{\YT}{$(Y,T)$\ }

\newcommand{\XtoYT}{$(X,T)\overset\pi\to(Y,T)$\ }

\newcommand{\OT}{\mathcal{O}_T}
\newcommand{\OS}{\mathcal{O}_S}
\newcommand{\OCT}{{\bar{\mathcal{O}}}_T}

\newcommand{\al}{\alpha}

\newcommand{\ga}{\gamma}
\newcommand{\del}{\delta}
\newcommand{\Del}{\Delta}
\newcommand{\ep}{\epsilon}
\newcommand{\sig}{\sigma}

\newcommand{\la}{\lambda}
\newcommand{\La}{\Lambda}
\newcommand{\tet}{\theta}

\newcommand{\om}{\omega}
\newcommand{\Om}{\Omega}

\newcommand{\ol}{\overline}

\newcommand{\br}{\vspace{3 mm}}
\newcommand{\imp}{\Rightarrow}

\newcommand{\wdis}{\curlywedge}

\newcommand{\nek}{,\ldots,}
\newcommand{\prodl}{\prod\limits}

\newcommand{\htop}{h_{{\rm top}}}

\newcommand{\inte}{{\rm{int\,}}}
\newcommand{\cls}{{\rm{cls\,}}}

\newcommand{\Aut}{{\rm{Aut\,}}}

\newcommand{\gr}{{\rm{gr\,}}}

\newcommand{\id}{{\rm{id}}}

\newcommand{\supp}{{\rm{supp\,}}}

\newcommand{\card}{{\rm{card\,}}}

\newcommand{\erg}{{\rm{erg}}}

\newcommand{\TE}{\mathbf{TE}}
\newcommand{\MIN}{\mathbf{MIN}}
\newcommand{\WM}{\mathbf{WM}}
\newcommand{\MM}{\mathbf{MM}}

\swapnumbers
\theoremstyle{plain}
\newtheorem{thm}{Theorem}[section]
\newtheorem{cor}[thm]{Corollary}
\newtheorem{lem}[thm]{Lemma}
\newtheorem{prop}[thm]{Proposition}

\theoremstyle{definition}
\newtheorem{defn}[thm]{Definition}

\newtheorem{rem}[thm]{Remark}

\begin{document}


\title[measurable and topological dynamics]
{On the interplay between measurable
and topological dynamics}

\author{E. Glasner and B. Weiss}

\address{Department of Mathematics\\
     Tel Aviv University\\
         Tel Aviv\\
         Israel}
\email{glasner@math.tau.ac.il}
\address {Institute of Mathematics\\
 Hebrew University of Jerusalem\\
Jerusalem\\
 Israel}
\email{weiss@math.huji.ac.il}

\begin{date}
{December 7, 2003}
\end{date}

\maketitle


\tableofcontents
\setcounter{secnumdepth}{2}

\newpage

\setcounter{section}{0}


\section*{Introduction}
Recurrent - wandering, conservative - dissipative, contracting -
expanding, deterministic - chaotic, isometric - mixing, periodic -
turbulent, distal - proximal, the list can go on and on. These
(pairs of) words --- all of which can be found in the dictionary
--- convey dynamical images and were therefore adopted by
mathematicians to denote one or another mathematical aspect of a
dynamical system.

The two sister branches of the theory of dynamical  systems called
{\em ergodic theory} (or {\em measurable dynamics}) and {\em
topological dynamics} use these words to describe different but
parallel notions in their respective theories and the surprising
fact is that many of the corresponding results are rather similar.
In the following  article we have tried to demonstrate both the
parallelism and the discord between ergodic theory and topological
dynamics. We hope that the subjects we chose to deal with will
successfully demonstrate this duality.

The table of contents gives a detailed listing of the topics
covered. In the first part we have detailed the strong analogies
between ergodic theory and topological dynamics as shown in the
treatment of recurrence phenomena, equicontinuity and weak
mixing, distality and entropy. In the case of distality the
topological version came first and the theory of measurable
distality was strongly influenced by the topological results.
For entropy theory the influence clearly was in the opposite
direction. The prototypical result of the second part is the
statement that any abstract measure probability preserving
system can be represented as a continuous transformation of a
compact space, and thus in some sense ergodic theory embeds into
topological dynamics.

\br

We have not attempted in any way to be either systematic or
comprehensive. Rather our choice of subjects was motivated by
taste, interest and knowledge and to great extent is random.
We did try to make the survey accessible to non-specialists,
and for this reason we deal throughout with the simplest case
of actions of $\Z$.
Most of the discussion carries over to noninvertible mappings
and to $\R$ actions. Indeed much of what we describe can be
carried over to general amenable groups.
Similarly, we have for the most part given rather complete definitions.
Nonetheless, we did take advantage of the fact that this article
is part of a handbook and for some of the definitions,
basic notions and well known results we refer the reader to the
earlier introductory chapters of volume I. Finally,
we should acknowledge the fact that  we made use  of parts of
our previous expositions
\cite{W4} and \cite{G}.

\br

We  made the writing of this survey more pleasurable for us by the
introduction of a few original results. In particular the
following results are entirely or partially new. Theorem
\ref{periodic} (the equivalence of the existence of a Borel
cross-section with the coincidence of recurrence and periodicity),
most of the material in Section 4 (on topological mild-mixing),
all of subsection 7.4 (the converse side of the local variational
principle) and subsection 7.6 (on topological determinism).

\br


{\large{\part{Analogies}}}

\br

\section{Poincar\'e recurrence vs. Birkhoff's
recurrence}\label{Sec-Poin}

\subsection{Poincar\'e recurrence theorem and
topological recurrence}
The simplest dynamical systems are the periodic ones.
In the absence of periodicity the crudest approximation
to this is approximate periodicity where instead of
some iterate $T^nx$ returning exactly to $x$ it returns
to a neighborhood of $x$. The first theorem in abstract
measure dynamics is Poincar\'{e}'s recurrence theorem
which asserts that for a finite measure preserving
system $(X,\Bcal,\mu,T)$ and any measurable set $A$,
$\mu$-a.e. point of $A$ returns to $A$
(see \cite[Theorem 4.3.1]{HKat}).
The proof of this
basic fact is rather simple and depends on identifying
the set of points $W\subset A$ that never return to $A$.
These are called the {\bf wandering points\/} and their
measurability follows from the formula
$$
W= A \cap \left( \bigcap_{k=1}^\infty T^{-k}
(X\setminus A)\right).
$$
Now for $n\ge 0$, the sets $T^{-n} W$ are pairwise
disjoint since $x\in T^{-n}W$ means that the forward
orbit of $x$ visits $A$ for the last time at moment $n$.
Since $\mu(T^{-n}W)=\mu(W)$ it follows that $\mu(W)=0$
which is the assertion of Poincar\'{e}'s theorem.
Noting that $A\cap T^{-n}W$ describes the points of $A$
which visit $A$ for the last time at moment $n$, and that
$\mu(\cup_{n=0}^\infty T^{-n} W)=0$ we have established the
following stronger formulation of Poincar\'{e}'s theorem.

\begin{thm}
For a finite measure preserving system $(X,\Bcal,\mu,T)$
and any measurable set $A$,
$\mu$-a.e. point of $A$ returns to $A$ infinitely often.
\end{thm}

\br

Note that only sets of the form $T^{-n}B$ appeared in the
above discussion so that the invertibility of $T$ is not
needed for this result. In the situation of classical dynamics,
which was Poincar\'{e}'s main interest, $X$ is also equipped
with a separable metric topology. In such a situation we can
apply the theorem to a refining sequence of partitions $\Pcal_m$,
where each $\Pcal_m$ is a countable partition into sets of
diameter at most $\frac{1}{m}$. Applying the theorem to a fixed
$\Pcal_m$ we see that $\mu$-a.e. point comes to within
$\frac{1}{m}$ of itself, and since the intersection of a
sequence of sets of full measure has full measure, we deduce
the corollary that $\mu$-a.e. point of $X$ is recurrent.

This is the measure theoretical path to the recurrence
phenomenon which depends on the presence of a finite invariant
measure. The necessity of such measure is clear from considering
translation by one on the integers. The system is dissipative,
in the sense that no recurrence takes place even though there is
an infinite invariant measure.

\br

{\begin{center}{$\divideontimes$}\end{center}}

\br

There is also a topological path to recurrence which was developed
in an abstract setting by G. D. Birkhoff. Here the above example
is eliminated by requiring that the topological space $X$,
on which our continuous transformation $T$ acts, be compact.
It is possible to show that in this setting a finite $T$-invariant
measure always exists, and so we can retrieve the measure
theoretical picture, but a purely topological discussion will
give us better insight.

\br

A key notion here is that of minimality. A nonempty closed,
$T$-invariant set $E\subset X$, is said to be {\bf minimal\/}
if $F\subset E$, closed and $T$-invariant implies
$F=\emptyset$ or $F=E$. If $X$ itself is a minimal set we
say that the system \XT is a {\bf minimal system}.

Fix now a point $x_0\in X$ and consider
$$
\om(x_0)=\bigcap_{n=1}^\infty\ol{\{T^k x_0:
k\ge n\}}.
$$
The points of $\om(x_0)$ are called {\bf $\om$-limit
points of $x_0$\/}, ($\om =$ last letter of the Greek
alphabet) and in the separable case $y\in \om(x_0)$ if
and only if there is some sequence $k_i\to \infty$ such that
$T^{k_i}x_0 \to y$. If $x_0\in \om(x_0)$ then $x_0$ is called
a {\bf positively recurrent point}.

Clearly $\om(x_0)$ is a closed and $T$-invariant set.
Therefore, in any nonempty minimal set $E$, any point
$x_0\in E$ satisfies $x_0\in \om(x_0)$ and thus we see
that minimal sets have recurrent points.

In order to see that compact systems \XT have recurrent points it
remains to show that minimal sets always exist. This is an immediate
consequence of Zorn's lemma applied to the family of nonempty
closed $T$-invariant subsets of $X$. A slightly more
constructive proof can be given when $X$ is a compact and
separable metric space. One can then list a sequence of
open sets $U_1, U_2, \dots$ which generate the topology,
and perform the following algorithm:

\begin{enumerate}
\item
set $X_0=X$,
\item
for $i=1,2, \dots$, \newline
if $\bigcup_{\,n=-\infty}^{\,\infty}
T^{-n}U_i\supset X_{i-1}$
put $X_i= X_{i-1}$, else put
$X_i= X_{i-1}\setminus \bigcup_{\,n=-\infty}^{\,\infty} T^{-n}U_i$.
\end{enumerate}

Note that $X_i\ne\emptyset$ and closed and thus
$X_\infty=\bigcap_{\,i=0}^{\,\infty} X_i$ is nonempty.
It is clearly $T$-invariant and for any $U_i$, if $U_i\cap X_\infty
\ne\emptyset$ then
$\bigcup_{\,-\infty}^{\,\infty} T^{-n}(U_i\cap X_\infty)
=X_\infty$, which shows that $(X_\infty, T)$ is minimal.

\br

\subsection{The existence of Borel cross-sections}

There is a deep connection between recurrent points in
the topological context and ergodic theory. To see this
we must consider quasi-invariant measures.
For these matters it is better to enlarge the scope and
deal with continuous actions of $\Z$, generated by $T$,
on a {\em complete separable metric space\/} $X$.
A probability measure $\mu$ defined on the Borel subsets
of $X$ is said to be {\bf quasi-invariant\/}  if $T\cdot \mu
\sim \mu$. Define such a system $(X,\Bcal,\mu,T)$ to be
{\bf conservative\/} if for any measurable set $A$,
$\ TA\subset A$ implies $\mu(A\setminus TA)=0$.

It is not hard to see that the conclusion of Poincar\'{e}'s
recurrence theorem holds for such systems; i.e.
if $\mu(A)>0$, then $\mu$-a.e. $x$ returns
to $A$ infinitely often. Thus once again $\mu$-a.e. point
is topologically recurrent. It turns out now that the
existence of a single topologically recurrent point
implies the existence of a non-atomic conservative
quasi-invariant measure. A simple proof of this fact
can be found in \cite{KW} for the case when $X$ is compact
--- but the proof given there is equally valid for complete
separable metric spaces. In this sense the phenomenon of
Poincar\'{e} recurrence and topological recurrence are
``equivalent" with each implying the other.

\br

A Borel set $B\subset X$ such that each orbit
intersects $B$ in exactly one point is called
{\bf a Borel cross-section\/} for the system \XT.
If a Borel cross-section exists, then no
non-atomic conservative quasi-invariant measure can
exist. In \cite{Wei} it is shown that the converse is also
valid --- namely if there are no conservative quasi-invariant
measures then there is a Borel cross-section.

Note that the periodic points of \XT form a Borel subset
for which a cross-section always exists, so that we can
conclude from the above discussion the following statement
in which no explicit mention is made of measures.

\begin{thm}\label{periodic}
For a system \XT, with $X$ a completely metrizable
separable space, there exists a Borel cross-section
if and only if the only recurrent points are the
periodic ones.
\end{thm}

\begin{rem}
Already in \cite{Glimm} as well as in \cite{Eff} one finds
many equivalent conditions for the existence of a Borel
section for a system \XT. However one doesn't find there explicit
mention of conditions in terms of recurrence.
Silvestrov and Tomiyama \cite{ST} established the theorem
in this formulation for $X$ compact (using $C^*$-algebra methods).
We thank A. Lazar for drawing our attention to their paper.
\end{rem}

\br

\subsection{Recurrence sequences and Poincar\'e sequences}
We will conclude this section with a discussion of
recurrence sequences and Poincar\'e sequences. First for
some definitions.
Let us say that $D$ is a {\bf recurrence set} if  for any
dynamical system $(Y, T)$ with compatible metric $\rho$ and
any $\epsilon >0$ there is a point $y_0$ and a $d \in D$ with
$$
\rho(T^dy_0,\ y_0)< \epsilon.
$$
Since any system contains minimal sets it suffices to restrict
attention here to minimal systems.  For minimal systems the set of
such $y$'s for a fixed $\epsilon$ is a dense open set.

To see this fact, let $U$ be an open set.  By the minimality there
is some $N$ such that for any $y \in Y$, and some $0 \leq n \leq
N$, we have $T^ny \in U$.  Using the uniform continuity of $T^n$,
we find now a $\delta >0$ such that if $\rho(u,\ v) < \delta$ then
for all $0 \leq n \leq N$
$$
\rho(T^nu,\ T^n v)< \epsilon.
$$
Now let $z_0$ be a point in $Y$ and $d_0 \in D$ such that
\begin{equation}\label{d}
p(T^{d_0}z_0,\ z_0) < \delta.
\end{equation}
For some $0 \leq n_0 \leq N$ we have $T^{n_0}z_0 = y_0 \in U$ and
from \eqref{d} we get $\rho(T^{d_0}y_0,\ y_0) < \epsilon$. Thus points
that $\epsilon$ return form an open dense set. Intersecting over
$\epsilon\to 0$ gives a dense $G_\delta$ in $Y$ of points $y$ for
which
$$
\inf_{d\in D} \ \rho(T^dy,\ y)=0.
$$
Thus there are points which actually recur along times drawn from
the given recurrence set.

A nice example of a recurrence set is the set of squares.
To see this it is easier to prove a stronger  property which is
the analogue in ergodic theory of recurrence sets.

\begin{defn}
A sequence $\{s_j\}$ is said to be a
{\bf Poincar\'e sequence\/} if for any finite measure preserving system
$(X,\ \Bcal ,\ \mu,\ T)$ and any $B \in \Bcal$ with
positive measure we have
$$
\mu ( T ^{s_j} B \cap B ) > 0 \qquad {\text{for  some $ s_j$
in the sequence.}}
$$
\end{defn}

Since any minimal topological system $(Y,T)$ has finite
invariant measures with global support, $\mu$ any Poincar\'e
sequence is recurrence sequence. Indeed for any presumptive
constant $b>0$ which would witness the non-recurrence of $\{s_j\}$
for $(Y,T)$, there would have to be an open  set $B$
with diameter less than $b$ and having positive $\mu$-measure such
that  $T ^{s_j} B \cap B$ is empty for all $\{s_j\}$.

Here is a sufficient condition for a sequence to be a Poincar\'e
sequence:

\begin{lem}
If for every $\alpha \in (0,\ 2 \pi)$
$$
\lim_{n \to \infty}\ \frac {1}{n}\ \sum^{n}_{k=1}\
e^{i \alpha s_k}\ =\ 0
$$
then $\{s_k\}_1^\infty$ is a Poincar\'e sequence.
\end{lem}

\begin{proof}
Let $(X,\ \Bcal ,\ \mu,\ T)$ be a measure preserving
system and let $U$ be the unitary operator defined on $L^2(X,\
\Bcal,\ \mu)$ by the action of $T$, i.e.
$$
(Uf)(x)= f(Tx).
$$
Let $H_0$ denote the subspace of invariant functions and
for a set of positive measure $B$, let $f_0$ be the projection of
$1_B$ on the invariant functions.  Since this can also be seen as
a conditional expectation with respect to the $\sigma$-algebra of
invariant sets $f_0 \geq 0$ and is not zero.  Now since $\ch_B -
f_0$ is orthogonal to the space of invariant functions its
spectral measure with respect to $U$ doesn't have any atoms at
$\{0\}$.  Thus from the spectral representation we deduce that in
$L^2$-norm
$$
\left| \left| \frac {1}{n}\ \sum_1^n\ U^{s_k}
(1_B-f_0)\right|\right|_{L^2} \longrightarrow 0
$$
or
$$
\left| \left| \left( \frac {1}{n} \sum^n_1 \ U ^{s_k}\ 1_B
\right)-f_0 \right|\right|_{L_2} \longrightarrow 0
$$
and integrating against $1_B$ and using the fact that $f_0$ is the
projection of $1_B$ we see that
$$
\lim_{n \to \infty}\ \frac {1}{n}\ \sum^n_1\ \mu(B \cap
T^{-s_k}B) = \|f_0\|^2 >0
$$
which clearly implies that $\{s_k\}$ is a Poincar\'e sequence.
\end{proof}

The proof we have just given is in fact von-Neumann's original
proof for the mean ergodic theorem.  He used the fact that $\mathbb
{N}$ satisfies the assumptions of the proposition, which is Weyl's
famous theorem on the equidistribution of $\{n \alpha\}$.
Returning to the squares Weyl also showed that $\{n^2 \alpha\}$ is
equidistributed for all irrational $\alpha$. For rational $\alpha$
the exponential sum in the lemma needn't vanish , however the
recurrence along squares for the rational part of the  spectrum is
easily verified directly so that we can conclude that indeed the
squares are a Poincar\'e sequence and hence a recurrence sequence.

The converse is not always true, i.e. there are recurrence
sequences that are not Poincar\'e sequences. This was first
shown by I. Kriz  \cite{Kr} in a beautiful example
(see also \cite[Chapter 5]{W4}). Finally here is a
simple problem.

\br

{\bf Problem:}\   If $D$ is a recurrence sequence for all
circle rotations is it a recurrence set?

\br

A little bit of evidence for a positive answer to that problem
comes from looking at a slightly different characterization of
recurrence sets.  Let $\Ncal $ denote the collection of sets of
the form
$$
N(U,\ U)=\{n:\ T^{-n}U \cap U \neq \emptyset\},\qquad (U\ {\text{open
and nonempty}}),
$$
where $T$ is a minimal transformation.  Denote by $\Ncal^*$ the
subsets of $\mathbb {N}$ that have a non-empty intersection with every
element of $\Ncal $.  Then $\Ncal^*$ is exactly the class of
recurrence sets.  For minimal transformations, another description
of $N(U,\ U)$ is obtained by fixing some $y_0$ and denoting
$$
N(y_0,\ U)=\{n:\ T^ny_0 \in U\}
$$
Then $N(U,\ U)=N(y_0,\ U)-N(y_0,\ U)$.  Notice that the minimality
of $T$ implies that $N(y_0,\ U)$ is a {\bf syndetic} set (a set
with bounded gaps) and so any $N(U,\ U)$ is the set of differences
of a syndetic set.  Thus $\Ncal$ consists essentially of all sets
of the form $S\ -\ S$ where $S$ is a syndetic set.

Given a finite set of real numbers
$\{\la_1,\la_2,\dots,\la_k\}$ and $\ep>0$
set
$$
V(\la_1,\la_2,\dots,\la_k;\ep)=
\{n\in \Z: \max_{j} \{\|n\la_j \|<\ep\}\},
$$
where $\|\cdot\|$ denotes the distance to the closest integer.
The collection of such sets forms a basis of neighborhoods at zero
for a topology on $\Z$ which makes it a topological group.
This topology is called the {\bf Bohr topology}.
(The corresponding uniform structure is
totally bounded and the completion of $\Z$
with respect to it is a compact topological group
called the {\bf Bohr compactification} of $\Z$.)

Veech proved in \cite{Veech} that any set
of the form $S\ -\ S$ with $S\subset \Z$ syndetic contains
a neighborhood of zero in the Bohr topology
{\em up to a set of zero density}.
It is not known if in that statement the zero density set can be
omitted.  If it could then a positive answer to the above problem
would follow (see also \cite{G26}).

\br

\section{The equivalence of weak mixing and continuous
spectrum}\label{Sec-WM}

In order to analyze the structure of a dynamical system
$\Xb$ there are, a priori, two possible approaches. In the first
approach one considers the collection of {\bf subsystems} $Y\subset X$
(i.e. closed $T$-invariant subsets)
and tries to understand how $X$ is built up by these subsystems.
In the other approach one is
interested in the collection of {\bf factors}
$X\overset{\pi}{\to}Y$ of the system $\Xb$.
In the measure theoretical case
the first approach leads to the ergodic decomposition and
thereby to the study of the ``indecomposable" or ergodic components
of the system. In the topological setup there is, unfortunately,
no such convenient decomposition describing the system in
terms of its indecomposable parts and one has to use some less
satisfactory substitutes. Natural candidates for indecomposable
components of a topological dynamical system are the ``orbit
closures" (i.e. the topologically transitive subsystems) or
the ``prolongation" cells (which often coincide with the
orbit closures), see \cite{AG0}. The minimal subsystems are of particular
importance here. Although we can not say, in any reasonable sense,
that the study of the general system can be reduced to that
of its minimal components, the analysis of the minimal systems
is nevertheless an important step towards a better understanding
of the general system.

This reasoning leads us to the study of the collection of
indecomposable systems (ergodic systems
in the measure category and transitive or minimal systems
in the topological case) and their factors. The simplest
and best understood indecomposable dynamical systems are the
ergodic translations of a compact monothetic group (a
cyclic permutation on $\Z_p$ for a prime number $p$,
the ``adding machine" on $\prod_{n=0}^\infty \Z_2$,
an irrational rotation $z\mapsto e^{2\pi i\al}z$ on
$S^1=\{z\in \C:|z|=1\}$ etc.). It is not hard to show that
this class of ergodic actions is characterized as those
dynamical systems which admit a model
$(X,\Xcal,\mu,T)$ where $X$ is a compact metric space,
$T:X\to X$ a surjective isometry and $\mu$ is $T$-ergodic.
We call these systems
{\bf Kronecker\/} or {\bf isometric\/}  systems.
Thus our first question
concerning the existence of factors should be: given
an ergodic dynamical system $\Xb$ which are its Kronecker factors?
Recall that a measure dynamical system $\Xb=(X,\Xcal,\mu,T)$
is called {\bf weakly mixing\/} if the product system
$(X\times X,\Xcal\otimes\Xcal,\mu\times \mu,T\times T)$
is ergodic. The following classical theorem is
due to von Neumann. The short and elegant proof we
give was suggested by Y. Katznelson.

\begin{thm}\label{easy}
An ergodic system $\Xb$ is weakly mixing iff it admits no
nontrivial Kronecker factor.
\end{thm}
\begin{proof}
Suppose $\Xb$ is weakly mixing and admits an isometric
factor. Now a factor of a weakly mixing
system is also weakly mixing and the only system which is both
isometric and weakly mixing is the trivial system (an easy
exercise).
Thus a weakly mixing system does not admit a nontrivial
Kronecker factor.

For the other direction,
if $\Xb$ is non-weakly mixing then in the product space $X\times X$
there exists a $T$-invariant measurable subset $W$
such that $0<(\mu\times \mu)(W)<1$. For every $x\in X$ let
$W(x)=\{x'\in X:(x,x')\in W\}$ and let $f_x={\ch}_{W(x)}$,
a function in $L^\infty (\mu)$.
It is easy to check that $U_T f_x=f_{T^{-1}x}$ so that the
map $\pi:X\to L^2(\mu)$ defined by $\pi(x)=f_x, x\in X$
is a Borel factor map. Denoting
$$
\pi(X)=Y\subset L^2(\mu), \quad {\text{and}} \quad \nu=\pi_*(\mu),
$$
we now have a factor map $\pi: \Xb \to (Y,\nu)$.
Now the function $\|\pi(x)\|$ is clearly measurable
and invariant and by ergodicity it is a constant
$\mu$-a.e.; say $\|\pi(x)\|=1$.
The dynamical system $(Y,\nu)$ is thus a subsystem
of the compact dynamical system $(B,U_T)$, where
$B$ is the unit ball of the Hilbert space $L^2(\mu)$
and $U_T$ is the Koopman unitary operator induced by $T$
on $L^2(\mu)$. Now it is well known (see e.g.
\cite{G}) that a compact topologically
transitive subsystem which carries an invariant probability
measure must be a Kronecker system and our proof is complete.
\end{proof}

\br

Concerning the terminology we used in the proof of
Theorem \ref{easy}, B. O. Koopman, a student of G. D. Birkhoff and a
co-author of both
Birkhoff and von Neumann introduced the crucial idea of
associating with a measure dynamical system
$\Xb=(X,\Xcal,\mu,T)$ the unitary operator
$U_T$ on the Hilbert space $L^2(\mu)$.
It is now an easy matter to
see that Theorem \ref{easy} can be re-formulated
as saying that the system $\Xb$ is weakly mixing iff
the point spectrum of the Koopman operator
$U_T$ comprises the single complex
number $1$ with multiplicity $1$.
Or, put otherwise, that the one dimensional space of constant
functions is the eigenspace corresponding to the eigenvalue
$1$ (this fact alone is equivalent to the ergodicity of the
dynamical system) and that the restriction of $U_T$
to the orthogonal complement of the space of constant functions
has a continuous spectrum.

\br

{\begin{center}{$\divideontimes$}\end{center}}

\br

We now consider a topological analogue of this theorem.
Recall that a topological system \XT is {\bf topologically
weakly mixing\/}
when the product system $(X\times X,T\times T)$ is topologically
transitive. It is {\bf equicontinuous\/} when the family
$\{T^n:n\in \Z\}$ is an equicontinuous family of maps.
Again an equivalent condition is the existence of a
compatible metric with respect to which $T$ is an isometry.
And, moreover, a minimal system is equicontinuous iff
it is a minimal translation on a compact monothetic group.
We will need the following lemma.

\begin{lem}\label{cont}
Let \XT be a minimal system and $f:X\to \R$ a
$T$-invariant function with at least one point of continuity
(for example this is the case when $f$ is lower or
upper semi-continuous or more generally when it is the
pointwise limit of a sequence of continuous functions),
then $f$ is a constant.
\end{lem}

\begin{proof}
Let $x_0$ be a continuity point and $x$ an arbitrary point
in $X$. Since $\{T^n x:n \in \Z\}$ is dense and as the value
$f(T^n x)$ does not depend on $n$ it follows that
$f(x)=f(x_0)$.
\end{proof}

\br

\begin{thm}\label{eswdisp}
Let \XT be a minimal system
then \XT is topologically weakly mixing iff it
has no non-trivial equicontinuous factor.
\end{thm}
\begin{proof}
Suppose \XT is minimal and topologically weakly mixing and let
$\pi:(X,T)\to(Y,T)$ be an equicontinuous factor.
If $(x,x')$ is a point whose $T\times T$ orbit
is dense in $X\times X$ then $(y,y')=(\pi(x),\pi(x'))$
has a dense orbit in $Y\times Y$. However, if \YT is
equicontinuous then $Y$ admits a compatible metric
with respect to which $T$ is an isometry and the
existence of a transitive point in $Y\times Y$
implies that $Y$ is a trivial one point space.

Conversely,
assuming that $(X\times X,T\times T)$ is not transitive we will
construct an equicontinuous factor $(Z,T)$ of $(X,T)$.
As $(X,T)$ is a minimal system, there exists a $T$-invariant
probability measure $\mu$ on $X$ with full support.
By assumption there exists an open $T$-invariant subset $U$
of $X\times X$, such  that $\cls U:= M \subsetneq X\times X$.
By minimality the projections of $M$ to both $X$ coordinates
are onto.
For every $y\in X$ let $M(y)=\{x\in X:(x,y)\in M\}$,
and let $f_y=\ch_{M(y)}$ be the indicator function of the set
$M(y)$, considered as an element of $L^1(X,\mu)$.

Denote by $\pi:X\to L^1(X,\mu)$ the map $y\mapsto f_y$.
We will show that $\pi$ is a continuous
homomorphism, where we consider $L^1(X,\mu)$ as a dynamical
system with the isometric action of the group
$\{U^n_T:n\in \Z\}$ and $U_Tf(x)=f(Tx)$.
Fix $y_0\in X$ and $\epsilon>0$.
There exists an open neighborhood $V$ of the closed set
$M(y_0)$ with $\mu(V\setminus M(y_0))<\epsilon$.
Since $M$ is closed the set map $y\mapsto M(y), X\to 2^X$ is
upper semi-continuous and we can find a neighborhood $W$ of
$y_0$ such that $M(y)\subset V$ for every $y\in W$.
Thus for every $y\in W$ we have $\mu(M(y)\setminus M(y_0))<\epsilon$.
In particular,
$\mu(M(y))\le\mu(M(y_0))+\epsilon$ and it follows that the map
$y\mapsto \mu(M(y))$ is upper semi-continuous.
A simple computation shows that
it is $T$-invariant, hence, by Lemma \ref{cont}, a constant.

With $y_0,\epsilon$ and $V, W$ as above,
for every $y\in W$, $\mu(M(y)\setminus M(y_0))<\epsilon$ and
$\mu(M(y))=\mu(M(y_0))$, thus
$\mu(M(y)\Delta M(y_0))<2\epsilon$,
i.e., $\| f_y-f_{y_0} \|_1<2\epsilon$.
This proves the claim that $\pi$ is continuous.

Let $Z=\pi(X)$ be the image of $X$ in $L^1(\mu)$. Since
$\pi$ is continuous, $Z$ is compact.
It is easy to see that the $T$-invariance of $M$ implies that
for every $n\in \Z$ and $y\in X$,\
$f_{T^{-n}y}=f_y\circ T^n$ so that
$Z$ is $U_T$-invariant and
$\pi:(Y,T)\to (Z,U_T)$ is a homomorphism.
Clearly $(Z,U_T)$ is minimal and equicontinuous (in fact
isometric).
\end{proof}

\br

Theorem \ref{eswdisp} is due to
Keynes and Robertson \cite{KRo} who developed an idea of
Furstenberg, \cite{Fur2};
and independently to K. Petersen \cite{Pe1} who utilized a
previous work of W. A. Veech, \cite{Veech}.
The proof we presented is an elaboration of a work of
McMahon \cite{McM2}
due to Blanchard, Host and Maass, \cite{BHM}.
We take this opportunity to point out a curious phenomenon
which recurs again and again. Some problems in topological
dynamics --- like the one we just discussed --- whose formulation
is purely topological, can be solved using the
fact that a $\Z$ dynamical system always carries an invariant
probability measure, and then employing a machinery provided
by ergodic theory. In several cases this approach is the only
one presently known for solving the problem. In the present
case however purely topological proofs exist, e.g. the
Petersen-Veech proof is one such.

\br

\section{Disjointness: measure vs. topological}\label{Sec-disj}
In the ring of integers $\Z$ two integers $m$ and $n$ have no
common factor if whenever $k|m$ and $k|n$ then $k=\pm 1$.
They are disjoint if  $m\cdot n$ is the least common multiple of
$m$ and $n$.
Of course in $\Z$ these two notions coincide. In his seminal paper
of 1967 \cite{Fur3},
H. Furstenberg introduced the same notions in the
context of dynamical systems, both measure-preserving
transformations and homeomorphisms
of compact spaces, and asked whether in these categories as well
the two are equivalent.
The notion of a factor in, say the measure category, is the natural
one: the dynamical system
$\Yb=(Y,\Ycal,\nu,T)$ is a {\bf factor\/} of the dynamical system
$\Xb=(X,\Xcal,\mu,T)$ if there exists a measurable map
$\pi:X \to Y$ with $\pi(\mu)=\nu$ that
$T\circ \pi= \pi \circ T$.
A common factor of two systems $\Xb$ and $\Yb$ is thus
a third system $\Zb$ which is a factor of both.
A {\bf joining\/} of the two systems $\Xb$ and $\Yb$ is
any system $\Wb$ which admits both as factors and is in turn
spanned by them. According to Furstenberg's definition
the systems $\Xb$ and $\Yb$ are {\bf disjoint\/} if the product
system $\Xb\times \Yb$ is the only joining they admit.
In the topological category, a joining of \XT and $(Y,S)$
is any subsystem $W\subset X\times Y$ of the product system
$(X\times Y, T\times S)$ whose projections on both
coordinates are full; i.e. $\pi_X(W)=X$ and $\pi_Y(W)=Y$.
\XT and $(Y,S)$ are {\bf disjoint\/} if $X\times Y$
is the unique joining of these two systems.
It is easy to verify that if \XT and $(Y,S)$ are disjoint
then at least one of them is minimal. Also, if both systems
are minimal then they are disjoint iff the product
system $(X\times Y, T\times S)$ is minimal.

In 1979, D. Rudolph, using joining techniques,
provided the first example of a pair of ergodic
measure preserving transformations with no common
factor which are not disjoint \cite{Ru1}.
In this work Rudolph laid the foundation
of joining theory. He introduced the class of
dynamical systems having ``minimal self-joinings"
(MSJ), and constructed a rank one mixing dynamical system
having minimal self-joinings of all orders.

Given a dynamical system $\Xb=(X,\Xcal,\mu,T)$
a probability measure $\la$ on the
product of $k$ copies of $X$ denoted $X_1,X_2\nek X_k$,
invariant under the product transformation and
projecting onto $\mu$ in each coordinate
is a {\bf $k$-fold self-joining\/}. It is called an
{\bf off-diagonal\/} \label{def-offdiag} if it is
a ``graph" measure of the form
$\la={\gr}(\mu,T^{n_1},\dots,T^{n_k})$,
i.e. $\la$ is the image of $\mu$ under the map
$x\mapsto\big(T^{n_1}x,T^{n_2}x\nek T^{n_k}x\big)$
of $X$ into $\prodl^k_{i=1}X_i$.
The joining $\la$ is a {\bf product of off-diagonals\/}
if there exists a
partition $(J_1\nek J_m)$ of $\{1\nek k\}$ such that
(i)\ For each $l$, the projection of $\la$ on
$\prodl_{i\in J_l}X_i$ is an off-diagonal,
(ii)\ The systems $\prodl_{i\in J_l}X_i$,
$1\le l\le m$, are independent.
An ergodic system $\Xb$ has
{\bf minimal self-joinings of order $k$\/}
if every $k$-fold ergodic self-joining of $\Xb$ is a
product of off-diagonals.

In \cite{Ru1} Rudolph shows how any dynamical system with MSJ
can be used to construct a counter example
to Furstenberg's question as well as
a wealth of other counter examples to various questions
in ergodic theory.
In \cite{JRS} del Junco, Rahe and Swanson were able
to show that the classical example
of Chac{\'o}n \cite{Chacon} has MSJ,
answering a question of Rudolph whether a weakly
but not strongly mixing system with MSJ exists.
In \cite{GW1}  Glasner and Weiss provide a topological
counterexample, which also serves as a natural
counterexample in the measure category.
The example consists of two
horocycle flows which have no nontrivial common factor but are
nevertheless not disjoint.
It is based on deep results of
Ratner \cite{Rat} which provide a complete
description of the self joinings
of a horocycle flow.
More recently an even more striking example was given in the
topological category by E. Lindenstrauss, where two minimal
dynamical systems with no nontrivial factor share a common almost
1-1 extension, \cite{Lis1}.

Beginning with the pioneering works of Furstenberg and Rudolph,
the notion of joinings was exploited by many authors; Furstenberg
1977 \cite{Fur4}, Rudolph 1979 \cite{Ru1}, Veech 1982 \cite{V2},
Ratner 1983 \cite{Rat}, del Junco and Rudolph 1987 \cite{JR1},
Host 1991 \cite{Host},
King 1992 \cite{King2}, Glasner, Host and Rudolph 1992 \cite{GHR},
Thouvenot 1993 \cite{Thou2}, Ryzhikov 1994 \cite{Ry},
Kammeyer and Rudolph 1995 (2002) \cite{KR},
del Junco, Lema\'nczyk and Mentzen 1995 \cite{JLM},
and Lema\'nczyk, Parreau and Thouvenot 2000 \cite{LPT}, to mention
a few.
The negative answer to Furstenberg's question and the
consequent works on joinings and disjointness
show that in order to study the relationship between two
dynamical systems it is necessary to know all the possible joinings
of the two systems and to understand the nature of these joinings.

Some of the best known disjointness relations among families
of dynamical systems are the following:
\begin{itemize}
\item ${\id} \ \bot \  $ ergodic,
\item distal $\  \bot \ $ weakly mixing (\cite{Fur3}),
\item rigid $\  \bot \ $ mild mixing (\cite{FW2}),
\item zero entropy $\  \bot \ $ $K$-systems (\cite{Fur3}),
\end{itemize}
in the measure category and
\begin{itemize}
\item $F$-systems $ \ \bot \  $ minimal (\cite{Fur3}),
\item minimal distal $\  \bot \ $ weakly mixing,
\item minimal zero entropy $\  \bot \ $ minimal UPE-systems
(\cite{Bla2}),
\end{itemize}
in the topological category.

\br

\section{Mild mixing: measure vs. topological}\label{Sec-mm}

\begin{defn}\label{def-mild}
Let $\Xb=(X,\Xcal,\mu,T)$ be a measure dynamical system.
\begin{enumerate}
\item
The system $\Xb$ is {\bf rigid\/} if there exists a
sequence $n_k\nearrow \infty$  such that
$$
\lim \mu
\left(T^{n_k}A\cap A\right)=\mu (A)
$$
for every measurable subset
$A$ of $X$. We say that $\Xb$ is $\{n_k\}$-{\bf rigid\/}.
\item
An ergodic system is {\bf mildly mixing\/}
if it has no non-trivial rigid factor.
\end{enumerate}
\end{defn}

These notions were introduced in \cite{FW2}.
The authors show that the mild mixing property is
equivalent to the following multiplier property.

\begin{thm}
An ergodic system $\Xb=(X,\Xcal,\mu,T)$ is mildly mixing
iff for every ergodic (finite or infinite) measure preserving
system $(Y,\Ycal,\nu,T)$, the product system
$$
(X\times Y, \mu \times \nu, T\times T),
$$
is ergodic.
\end{thm}

Since every Kronecker system is rigid it follows from Theorem
\ref{easy} that mild mixing implies weak mixing.
Clearly strong mixing implies mild mixing.
It is not hard to construct rigid weakly mixing
systems, so that the class of mildly mixing systems
is properly contained in the class of weakly mixing
systems.
Finally there are mildly but not strongly mixing
systems; e.g. Chac{\'o}n's system is an example
(see Aaronson and Weiss \cite{AW}).

\br

We also have the following analytic characterization of
mild mixing.

\begin{prop}\label{ex-rigid-matrix}
An ergodic system $\Xb$ is mildly
mixing iff
$$
\limsup_{n\to\infty}\phi_f(n)<1,
$$
for every matrix coefficient $\phi_f$, where
for $f\in L^2(X,\mu), \|f\|=1$,\
$\phi_f(n):=
\langle U_{T^n} f,f\rangle$.
\end{prop}

\begin{proof}
If $\Xb\to \Yb$ is a rigid factor, then there exists
a sequence $ n_i\to \infty$ such that
$U_{T^{n_i}}\to {\id}$ strongly on $L^2(Y,\nu)$.
For any function $f\in L^2_0(Y,\nu)$ with $\|f\|=1$,
we have $\lim_{i\to\infty}\phi_f(n_i)=1$. Conversely, if
$\lim_{i\to\infty}\phi_f(n_i)=1$ for some
$n_i\nearrow \infty$ and $f\in L^2_0(X,\mu),
\|f\|=1$, then $\lim_{i\to\infty} U_{T^{n_i}}f=f$. Clearly
$f$ can be replaced by a bounded function and we let
$A$ be the sub-algebra of $L^\infty(X,\mu)$ generated
by $\{U_{T^n}f:n \in \Z \}$. The algebra
$A$ defines a non-trivial factor $\Xb\to \Yb$ such that
$U_{T^{n_i}}\to {\id}$ strongly on $L^2(Y,\nu)$.
\end{proof}

\br

We say that a collection $\Fcal$ of nonempty subsets of $\Z$ is a
{\bf family\/} if it is hereditary upward and {\bf proper\/} (i.e.
$A\subset B$ and $A\in \Fcal$ implies $B\in \Fcal$, and
$\Fcal$ is neither empty nor all of $2^\Z$).

With a family $\Fcal$ of nonempty subsets of $\Z$
we associate the {\bf dual family\/}
$$
\Fcal^{*}=\{E:E\cap F\ne\emptyset, \forall \ F\in \Fcal\}.
$$
It is easily verified that $\Fcal^{*}$ is indeed a family.
Also, for families, $\Fcal_1\subset \Fcal_2\
\imp\ \Fcal^{*}_1\supset \Fcal^{*}_2$, and $\Fcal^{**}=\Fcal$.

We say that a subset $J$ of $\Z$ has {\bf uniform
density 1\/} if for every $0<\la<1$ there exists
an $N$ such that for every interval $I\subset \Z$
of length $> N$ we have $|J \cap I| \ge \la |I|$.
We denote by $\Dcal$ the family of subsets of $\Z$
of uniform density 1.
It is also easy to see that $\Dcal$ has the
finite intersection property.

Let $\Fcal$ be a family of nonempty subsets of $\Z$
which is closed under finite intersections
(i.e. $\Fcal$ is a filter). Following
\cite{Fur5} we say that a sequence $\{x_n: n\in \Z\}$
in a topological space $X$ {\bf $\Fcal$-converges} to a point
$x\in X$ if for every neighborhood $V$ of $x$ the
set $\{n: x_n\in V\}$ is in $\Fcal$. We denote this by
$$
\Fcal\,{\text -}\,\lim x_n =x.
$$

We have the following characterization
of weak mixing for measure preserving systems
which explains more clearly its name.
\begin{thm}
The dynamical system $\Xb=(X,\Xcal,\mu,T)$ is
weakly mixing iff
for every $A, B \in \Xcal$ we have
$$
\Dcal\,{\text -}\,\lim \mu(T^{-n}A\cap B)
=\mu(A)\mu(B).
$$
\end{thm}

An analogous characterization of measure theoretical
mild mixing is obtained by considering the families of $I\!P$
and $I\!P^*$ sets.
An $I\!P$-{\bf set\/} is any subset of $\Z$ containing a subset
of the form
$I\!P\{n_i\}= \{n_{i_1}+n_{i_2}+\cdots+n_{i_k}:i_1<i_2<\cdots<i_k\}$,
for some infinite sequence $\{n_i\}_{i=1}^\infty$.
We let $\Ical$ denote the family of $I\!P$-sets and call
the elements of the dual family $\Ical^*$, {\bf $I\!P^*$-sets\/}.
Again it is not hard to see that the family of $I\!P^*$-sets
is closed under finite intersections.
For a proof of the next theorem we refer to \cite{Fur5}.

\begin{thm}
The dynamical system $\Xb=(X,\Xcal,\mu,T)$ is
mildly mixing iff
for every $A, B \in \Xcal$ we have
$$
\Ical^*\,{\text -}\,\lim\mu(T^{-n}A\cap B)
=\mu(A)\mu(B).
$$
\end{thm}

\br

{\begin{center}{$\divideontimes$}\end{center}}

\br

We now turn to the topological category.
Let \XT be a topological dynamical system.
For two non-empty open sets $U,V\subset X$ and a point
$x\in X$ set
\begin{gather*}
N(U,V)=\{n\in \Z: T^n U\cap V\ne\emptyset\},\quad
N_+(U,V)=N(U,V)\cap \Z_+\\
{\text{and}} \qquad N(x,V)=\{n\in \Z: T^n x\in V\}.
\end{gather*}
Notice that sets of the form $N(U,U)$ are symmetric.

We say that \XT is {\bf topologically transitive\/}
(or just {\bf transitive\/}) if
$N(U,V)$ is nonempty whenever $U,V \subset X$ are two non-empty
open sets.
Using Baire's category theorem it is easy to see that (for metrizable
$X$) a system \XT is topologically transitive iff there exists
a dense $G_\del$ subset $X_0\subset X$ such that $\OCT(x)=X$
for every $x\in X_0$.

We define the family $\Fcal_{\rm thick}$ of {\bf thick sets\/}
to be the collection of sets which contain arbitrary long
intervals. The dual family $\Fcal_{\rm synd}=\Fcal_{\rm thick}^*$
is the collection of {\bf syndetic sets\/} --- those sets $A\subset \Z$
such that for some positive integer $N$ the intersection
of $A$ with every interval of length $N$ is nonempty.

Given a family $\Fcal$ we say that a topological
dynamical system $(X,T)$
is {\bf $\Fcal$-recurrent\/} if $N(A,A)\in \Fcal$
for every nonempty open set $A\subset X$.
We say that a dynamical system is {\bf $\Fcal$-transitive\/}
if $N(A,B)\in \Fcal$
for every nonempty open sets $A, B \subset X$.
The class of $\Fcal$-transitive systems is denoted by $\Ecal_{\Fcal}$.
E.g. in this notation the class of {\bf topologically mixing systems\/}
is $\Ecal_{\rm cofinite}$, where we call a subset $A\subset \Z$
co-finite when $\Z\setminus A$ is a finite set.
We write simply $\Ecal=\Ecal_{\rm{infinite}}$
for the class of {\bf recurrent
transitive\/} dynamical systems.
It is not hard to see that
when $X$ has no isolated points \XT is topologically transitive
iff it is recurrent transitive.
From this we then deduce that a weakly mixing system is
necessarily recurrent transitive.


In a dynamical system \XT a point $x\in X$  is a
{\bf wandering point\/} if there exists an open neighborhood
$U$ of $x$ such that the collection $\{T^n U: n\in \Z\}$ is
pairwise disjoint.

\begin{prop}\label{rec-tra}
Let $\XT$ be a topologically transitive dynamical system; then
the following conditions are equivalent:
\begin{enumerate}
\item
$\XT\in \Ecal_{\text{\rm infinite}}$.
\item
The recurrent points are dense in $X$.
\item
$\XT$ has no wandering points.
\item
The dynamical system $(X_\infty,T)$, the one point compactification
of the integers with translation and a fixed point at infinity, is not a
factor of $\XT$.
\end{enumerate}
\end{prop}

\begin{proof}
1 $\imp$ 4\
If $\pi:X\to X_\infty$ is a factor
map then, clearly $N(\pi^{-1}(0),\pi^{-1}(0))=\{0\}$.

4 $\imp$ 3\
If $U$ is a nonempty open wandering subset of
$X$ then $\{T^jU : j\in \Z\}\cup (X\setminus \bigcup
\{T^jU : j\in \Z\})$
is a partition of $X$. It is easy to see that this
partition defines a factor map $\pi : X\to X_\infty$.

3 $\imp$ 2\
This implication is a consequence of the following:

\begin{lem}
If the dynamical system $\XT$ has no wandering points then
the recurrent points are dense in $X$.
\end{lem}

\begin{proof}
For every $\del>0$ put
$$
A_\del=\{x\in X:\exists j\not=0,\ d(T^jx,x)<\del\}.
$$
Clearly $A_\del$ is an open set and we claim that it is dense. In fact
given $x\in X$ and $\ep>0$ there exists $j\not=0$ with
$$
T^jB_\ep(x)\cap B_\ep(x)\not=\emptyset.
$$
If $y$ is a point in this intersection then $d(T^{-j}y,y)<2\ep$.
Thus for $\ep<\del/2$ we have $y\in A_\del$ and $d(x,y)<\ep$.
Now by Baire's theorem
$$
A=\bigcap_{k=1}^\infty A_{1/k}
$$
is a dense $G_\del$ subset of $X$ and each point in $A$ is
recurrent.
\end{proof}

2 $\imp$ 1\
Given $U,V$ nonempty open subsets of $X$ and $k\in N(U,V)$
let $U_0$ be the nonempty open subset $U_0=U\cap T^{-k}V$. Check
that $N(U_0,U_0)+k\subset N(U,V)$. By assumption $N(U_0,U_0)$ is
infinite and a fortiori so is $N(U,V)$. This completes the proof
of Proposition \ref{rec-tra}.
\end{proof}

\br

%
%


A well known characterization of the class
$\WM$ of topologically weakly mixing systems is due
to Furstenberg:
\begin{thm}
$\WM=\Ecal_{\text{\rm thick}}$.
\end{thm}

Following \cite{AG} we call the systems in
$\Ecal_{\text{\rm synd}}$ {\bf topologically ergodic\/}
and write $\TE$ for this class.
This is a rich class as we can see from the following
claim from \cite{GW}. Here $\MIN$ is the class of minimal
systems and $\Eb$ the class of $E$-systems; i.e.
those transitive dynamical systems \XT for which there exists
a probability invariant measure with full support.

\begin{thm}\label{MIN<TE}
$\MIN, \Eb \subset \TE$.
\end{thm}

\begin{proof}
1.\
The claim for $\MIN$ is immediate by the well
known characterization of minimal systems:
\XT is minimal iff $N(x,U)$ is syndetic for
every $x\in X$ and nonempty open $U\subset X$.

2.\
Given two non-empty open sets $U,V$
in $X$, choose $k\in \Z$ with $T^kU\cap V\not=\emptyset$.
Next set $U_0=T^{-k}V\cap U$, and observe that $k+N(U_0,U_0)
\subset N(U,V)$. Thus it is enough to show that $N(U,U)$ is
syndetic for every non-empty open $U$.
We have to show that
$N(U,U)$ meets every thick subset $B\subset\Z$.
By Poincar\'e's recurrence theorem,
$N(U,U)$ meets every set of the form $A-A=\{n-m:n,m\in A\}$
with $A$ infinite.
It is an easy exercise to show that every thick set $B$ contains
some $D^+(A)=\{a_n-a_m:n>m\}$ for an infinite sequence
$A=\{a_n\}$. Thus $\emptyset\not= N(U,U)\cap \pm D^+(A)
\subset N(U,U)\cap \pm B$. Since $N(U,U)$ is symmetric,
this completes the proof.
\end{proof}

\br

We recall (see the previous section)
that two dynamical systems \XT and \YT are disjoint
if every closed $T\times T$-invariant subset of $X\times Y$
whose projections on $X$ and $Y$ are full, is necessarily the
entire space $X\times Y$. It follows easily that when
\XT and \YT are disjoint, at least one of them must be minimal.
If both \XT and \YT are minimal then they are disjoint iff
the product system is minimal. We say that \XT and \YT are
{\bf weakly disjoint\/} when the product system
$(X\times Y, T\times T)$ is transitive.
This is indeed a very
weak sense of disjointness as there are systems
which are weakly disjoint from themselves. In fact, by definition
a dynamical system is topologically weakly mixing iff it is weakly
disjoint from itself.

If $\Pb$ is a class of recurrent transitive dynamical systems
we let $\Pb^{\wdis}$ be the class of recurrent transitive
dynamical systems which are weakly disjoint from every member
of $\Pb$
$$
\Pb^{\wdis}=\{(X,T) : X\times Y\in \Ecal \ {\text
{\rm for every \ }} (Y,T)\in \Pcal \}.
$$
We clearly have $\Pb\subset \Qb\ \imp \Pb^{\wdis}\supset
\Qb^{\wdis}$ and $\Pb^{\wdis \wdis \wdis}= \Pb^{\wdis}$.

\br

For the discussion of topologically mildly mixing systems it will
be convenient to deal with families of subsets
of $\Z_+$ rather than $\Z$.
If $\Fcal$ is such a family then
$$
\Ecal_\Fcal=
\{(X,T) : N_+(A,B) \in \Fcal\
{\text {\rm for every nonempty open}}\ A,B \subset X\}.
$$
Let us call a subset of $\Z_+$ a $SI\!P$-{\bf set\/}
(symmetric $I\!P$-set), if it contains a subset of the form
$$
SI\!P\{n_i\}=\{n_\al-n_\beta > 0: \ n_\al,n_\beta\in I\!P\{n_i\}
\cup \{0\}\},
$$
for an $I\!P$ sequence $I\!P\{n_i\}\subset \Z_+$.
Denote by $\Scal$ the family of $SI\!P$ sets.
It is not hard to show that
$$
\Fcal_{\rm thick}\subset\Scal\subset\Ical,
$$
(see \cite{Fur5}).
Hence
$\Fcal_{\rm syndetic} \supset \Scal^* \supset \Ical^*$,
hence
$\Ecal_{{\text {\rm synd}}}\supset
\Ecal_{\Scal^*}\supset
\Ecal_{\Ical^*}$,
and finally
$$
\Ecal^{\wdis}_{{\text {\rm synd}}}\subset
\Ecal^{\wdis}_{\Scal^*}\subset
\Ecal^{\wdis}_{\Ical^*}.
$$

\br

\begin{defn}
A topological dynamical system
\XT is called {\bf topologically mildly mixing\/}
if it is in $\Ecal_{\Scal^*}$ and we denote the collection of
topologically mildly mixing systems by $\MM=\Ecal_{\Scal^*}$.
\end{defn}

\br

\begin{thm}\label{MM=E^}
A dynamical system is in $\Ecal$ iff it is weakly disjoint from
every topologically mildly mixing system:
$$
\Ecal=\MM^{\wdis}.
$$
And conversely it is topologically mildly mixing iff it is
weakly disjoint from every recurrent transitive system:
$$
\MM=\Ecal^{\wdis}.
$$
\end{thm}

\begin{proof}
1.\
Since $\Ecal_{\Scal^*}$ is nonvacuous (for example every
topologically mixing system is in $\Ecal_{\Scal^*}$), it follows
that every system in $\Ecal_{\Scal^*}^{\wdis}$ is in $\Ecal$.

Conversely, assume that  $\XT$ is in $\Ecal$
but $\XT \not\in \Ecal_{\Scal^*}^{\wdis}$, and we
will arrive at a contradiction.
By assumption there exists $\YT\in \Ecal_{\Scal^*}$ and a nondense
nonempty open invariant subset $W\subset X\times Y$.
Then $\pi_X(W)=O$ is a nonempty open invariant subset of $X$.
By assumption $O$ is dense in $X$.
Choose open nonempty sets $U_0\subset X$ and $V_0\subset Y$
with $U_0\times V_0\subset W$.
By Proposition \ref{rec-tra} there exists a recurrent point
$x_0$ in $U_0\subset O$. Then
there is a sequence $n_i\to\infty$ such that for the
$I\!P$-sequence $\{n_\al\}=I\!P\{n_i\}_{i=1}^\infty$,\
$I\!P{\text{-}}\lim T^{n_\al}x_0=x_0$ (see \cite{Fur5}).
Choose $i_0$ such that $T^{n_\al}x_0\in U_0$ for $n_\al\in
J=I\!P\{n_i\}_{i\ge i_0}$ and set $D=SI\!P(J)$.
Given $V$ a nonempty open subset of $Y$ we have:
$$
D\cap N(V_0,V) \not=\emptyset.
$$
Thus for some $\al,\beta$ and $v_0\in V_0$,
$$
T^{n_\al-n_\beta}(T^{n_\beta}x_0,v_0)=
(T^{n_\al}x_0,T^{n_\al-n_\beta}v_0)
\in (U_0\times V)\cap W.
$$
We conclude that
$$
\{x_0\}\times Y\subset \cls W.
$$

The fact that in an $\Ecal$ system the recurrent points are dense
together with the observation that $\{x_0\}\times Y\subset \cls W$
for every recurrent point $x_0\in O$, imply that
$W$ is dense in $X\times Y$, a contradiction.

\br

2.\
From part 1 of the proof we have $\Ecal=\Ecal_{\Scal^*}^{\wdis}$,
hence $\Ecal^{\wdis}= \Ecal_{\Scal^*}^{\wdis\wdis}
\supset\Ecal_{\Scal^*}$.

Suppose $\XT\in \Ecal$ but $\XT\not\in \Ecal_{\Scal^*}$,
we will show that $\XT\not\in \Ecal^{\wdis}$. There
exist $U,V\subset X$, nonempty open subsets and an $I\!P$-set
$I=I\!P\{n_i\}$ for a monotone increasing sequence $\{n_1<n_2<
\cdots\}$ with
$$
N(U,V)\cap D=\emptyset,
$$
where
$$
D=\{n_\al-n_\beta: n_\al,n_\beta\in I,\ n_\al>n_\beta\}.
$$
If $\XT$ is not topologically weakly mixing then
$X\times X\not\in \Ecal$
hence $\XT\not\in \Ecal^{\wdis}$. So we can assume that
$\XT$ is topologically weakly mixing. Now in $X\times X$
$$
N(U\times V,V\times U)= N(U,V)\cap N(V,U)=N(U,V)\cap -N(U,V),
$$
is disjoint from $D\cup -D$, and replacing $X$ by $X\times X$
we can assume that $N(U,V)\cap (D\cup -D)=\emptyset$.
In fact, if $X\in \Ecal^{\wdis}$ then $X\times Y \in \Ecal$
for every $Y \in \Ecal$, therefore
$X\times (X\times Y) \in \Ecal$ and we see that
also $X\times X \in \Ecal^{\wdis}$.

By going to a subsequence, we can assume that
$$
\lim_{k\to\infty} n_{k+1}-\sum_{i=1}^kn_i =\infty.
$$
in which case the representation
of each $n\in I$ as $n=n_\al=n_{i_1}+n_{i_2}+\cdots
+n_{i_k}; \ \al=\{i_1<i_2<\cdots<i_k\}$ is unique.

Next let $y_0\in\{0,1\}^\Z$ be the sequence $y_0=\ch_I$.
Let $Y$ be the orbit closure of $y_0$ in $\{0,1\}^\Z$
under the shift $T$, and let $[1]=\{y\in Y:y(0)=1\}$.
Observe that
$$
N(y_0,[1])=I.
$$
It is easy to check that
$$
I\!P{\text{-}}\lim T^{n_\al}y_0=y_0.
$$
Thus the system $(Y,T)$ is topologically transitive with
$y_0$ a recurrent point; i.e. $\YT\in \Ecal$.

We now observe that
$$
N([1],[1])= N(y_0,[1])- N(y_0,[1])=I-I= D\cup -D\cup \{0\}.
$$
If $X\times Y$ is topologically transitive then
in particular
\begin{gather*}
N(U\times [1],V\times [1])=N(U,V)\cap N([1],[1])
=\\
 N(U,V)\cap (D\cup -D\cup \{0\})=\ {\text{\rm infinite\ set}}.
\end{gather*}
But this contradicts our assumption.
Thus $X\times Y\not\in \Ecal$ and $\XT\not\in \Ecal^{\wdis}$.
This completes the proof.
\end{proof}

\br

We now have the following:

\begin{cor}
Every topologically mildly mixing system is weakly
mixing and topologically ergodic:
$$
\MM\subset \WM\cap \TE.
$$
\end{cor}

\begin{proof}
We have $\Ecal_{\Scal^*}\subset \Ecal
=\Ecal_{\Scal^*}^{\wdis}$,
hence for every $\XT\in \Ecal_{\Scal^*}$,
$X\times X \in \Ecal$ i.e. $\XT$ is topologically weakly mixing.
And, as we have already observed the inclusion
$\Fcal_{\rm syndetic} \supset \Scal^*$,
entails
$\TE=\Ecal_{{\text {\rm synd}}}\supset
\Ecal_{\Scal^*} = \MM$.
\end{proof}

\br

To complete the analogy with the measure theoretical
setup we next define a topological analogue of rigidity.
This is just one of several possible definitions
of topological rigidity and we refer to \cite{GM}
for a treatment of these notions.

\begin{defn}
A dynamical system
\XT is called {\bf uniformly rigid\/} if there exists a
sequence $n_k\nearrow \infty$  such that
$$
\lim_{k\to\infty}  \sup_{x\in X} d(T^{n_k}x,x)=0,
$$
i.e. $\lim_{k\to\infty} T^{n_k}=\id$
in the uniform topology on the group of homeomorphism
of $H(X)$ of $X$.
We denote by $\Rcal$ the collection of topologically
transitive uniformly rigid systems.
\end{defn}

\br

In \cite{GM} the existence of minimal weakly mixing
but nonetheless uniformly rigid dynamical systems
is demonstrated. However, we have the following:

\begin{lem}
A system which is both topologically mildly mixing and
uniformly rigid is trivial.
\end{lem}

\begin{proof}
Let $\XT$ be both topologically mildly mixing and uniformly rigid.
Then
$$
\La=\cls \{T^n: n\in \Z\}\subset H(X),
$$
is a Polish monothetic group.

Let $T^{n_i}$ be a sequence converging uniformly to $\id$, the
identity element of $\La$. For a subsequence we can ensure
that $\{n_\al\}=I\!P\{n_i\}$ is an $I\!P$-sequence such that
$I\!P{\text{-}}\lim T^{n_\al} = \id$ in $\La$.
If $X$ is nontrivial we can now find an open ball
$B=B_\del(x_0)\subset X$ with $TB\cap B=\emptyset$.
Put $U=B_{\del/2}(x_0)$ and $V=TU$; then by assumption
$N(U,V)$ is an $SI\!P^*$-set and in particular:
$$
\forall \al_0\  \exists \al,\beta>\al_0,\ n_\al-n_\beta\in N(U,V).
$$
However, since $I\!P{\text{-}}\lim T^{n_\al} = \id$,
we also have eventually,
$T^{n_\al-n_\beta}U\subset B$; a contradiction.
\end{proof}

\br

\begin{cor}
A topologically mildly mixing system has no nontrivial
uniformly rigid factors.
\end{cor}

\br

We conclude this section with the following result
which shows how these topological and measure theoretical
notions are related.

\begin{thm}
Let \XT be a topological dynamical system with the
property that there exists an invariant probability measure
$\mu$ with full support such that the
associated measure preserving dynamical system
$(X,\Xcal,\mu,T)$ is  measure theoretically
mildly mixing then \XT is topologically mildly mixing.
\end{thm}

\begin{proof}
Let $(Y,S)$ be any system in $\Ecal$; by
Theorem \ref{MM=E^} it suffices to show that
$(X\times Y,T\times S)$ is topologically transitive.
Suppose $W\subset
X\times Y$ is a closed $T\times S$-invariant set with
$\inte W\ne\emptyset$.
Let $U\subset X, V\subset V$ be two nonempty open
subsets with $U\times V\subset W$. By transitivity of $(Y,S)$
there exits a transitive recurrent point $y_0\in V$.
By theorems of Glimm and Effros
\cite{Glimm}, \cite{Eff}, and Katznelson and Weiss \cite{KW} (see also
Weiss \cite{Wei}),
there exists a (possibly infinite) invariant
ergodic measure $\nu$ on $Y$ with $\nu(V)>0$.

Let $\mu$ be the probability invariant measure of full
support on $X$ with respect to which $(X,\Xcal,\mu,T)$ is
measure theoretically mildly mixing.
Then by \cite{FW2} the measure $\mu\times \nu$ is ergodic.
Since $\mu\times\nu(W) \ge \mu\times\nu(U\times V) >0$
we conclude that $\mu\times\nu(W^c)=0$ which clearly implies
$W= X\times Y$.
\end{proof}

\br

We note that the definition of topological mild mixing
and the results described above concerning this notion
are new. However independently of our work
Huang and Ye in a recent work also define
a similar notion and give it a comprehensive
and systematic treatment, \cite{HY2}. The first named
author would like to thank E. Akin for instructive
conversations on this subject.

Regarding the classes $\WM$ and $\TE$ let us mention the
following result from \cite{W}.

\begin{thm}
$$
\TE = \WM^{\wdis}.
$$
\end{thm}

For more on these topics we refer to
\cite{Fur5}, \cite{A}, \cite{W}, \cite{AG}, \cite{HY1}
and \cite{HY2}.

\br

\section{Distal systems: topological vs. measure}\label{Sec-distal}
As noted above the Kronecker or minimal equicontinuous dynamical
systems can be considered as the most elementary type of systems.
What is then the next stage? The clue in the topological case,
which chronologically came first, is to be found in the notion
of distality.
A topological system \XT is called {\bf distal \/} if
$$
\inf_{n\in \Z} d(T^nx,T^nx')>0
$$
for every $x\ne x'$ in $X$.
It is easy to see that this property does not depend on
the choice of a metric. And, of course, every equicontinuous
system is distal. Is the converse true? Are these notions
one and the same? The dynamical system given on the unit
disc $D=\{z\in \C:|z|\le 1\}$ by the formula $Tz=z\exp(2\pi i |z|)$
is a counter example, it is distal but not equicontinuous.
However it is not minimal. H. Furstenberg in 1963 noted that
skew products over an equicontinuous basis with compact group
translations as fiber maps are always distal, often minimal,
but rarely equicontinuous, \cite{Fur2}.
A typical example is the homeomorphism
of the two torus $\T^2=\R^2/\Z^2$ given by $T(x,y)=(x+\al,y+x)$
where $\al\in \R/Z$ is irrational. Independently and at about the
same time, it was shown by L. Auslander, L. Green and F. Hahn
that minimal nilflows are distal
but not equicontinuous, \cite{AGH}.
These examples led Furstenberg to his path
breaking structure theorem, \cite{Fur2}.

Given a homomorphism $\pi:(X,T)\to (Y,T)$ let
$R_\pi=\{(x,x'):\pi(x)=\pi(x')\}$. We say that the homomorphism
$\pi$ is an {\bf isometric extension\/} if there exists a
continuous function $d:R_\pi\to \R$ such that for each
$y\in Y$ the restriction of $d$ to $\pi^{-1}(y)\times \pi^{-1}(y)$
is a metric and for every $x,x'\in \pi^{-1}(y)$ we have
$d(Tx,Tx')=d(x,x')$.

If $K$ is a compact subgroup of ${\Aut}(X,T)$ (the group
of homeomorphisms of $X$ commuting with $T$, endowed
with the topology of uniform convergence)
then the map $x\mapsto Kx$ defines a factor map
\XtoYT with $Y=X/K$ and
$R_\pi=\{(x,kx):x\in X,\ k\in K\}$.
Such an extension is called a
{\bf group extension\/}\label{def-gr-ext-t}.
It turns out, although this is not so easy
to see, that when \XT is minimal then $\pi:(X,T)\to (Y,T)$
is an isometric extension iff there exists a commutative diagram:
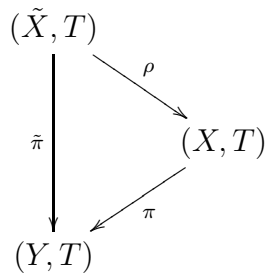
\begin{equation*}\label{iso-gr-diag}
\xymatrix
{
(\tilde X,T) \ar[dd]_{\tilde\pi} \ar[dr]^{\rho}  &  \\
& (X,T) \ar[dl]^{\pi}\\
(Y,T) &
}
\end{equation*}
where $(\tilde X,T)$ is minimal and
$(\tilde X,T)\overset{\tilde\pi}\to(X,T)$
is a group extension with some compact group
$K\subset \Aut(\tilde X,T)$
and the map $\rho$ is the quotient map from
$\tilde X$ onto $X$ defined by a closed subgroup $H$ of $K$.
Thus $Y=\tilde X/K$ and $X=\tilde X/H$ and we can think
of $\pi$ as a {\bf homogeneous space extension\/} with fiber
$K/H$.

We say that a (metrizable) minimal system \XT is an {\bf $I$ system}
if there is a (countable) ordinal $\eta$ and a family of systems
$\{(X_\tet,x_\tet)\}_{\tet\le\eta}$
such that (i) $X_0$ is the trivial system, (ii) for every
$\tet<\eta$ there exists an isometric homomorphism
$\phi_\tet:X_{\tet+1}\to X_\tet$,
(iii) for a limit ordinal $\la\le\eta$ the system $X_\la$
is the inverse limit of the systems $\{X_\tet\}_{\tet<\la}$
(i.e. $X_\la=\bigvee_{\tet<\la}(X_\tet,x_\tet)$),
and
(iv) $X_\eta=X$.

\begin{thm}[Furstenberg's structure theorem]
\label{furst-structure-tm}
A minimal system is distal iff it is an I-system.
\end{thm}

\br

{\begin{center}{$\divideontimes$}\end{center}}

\br

W. Parry in his 1967 paper \cite{Pa1} suggested an intrinsic
definition of measure distality. He defines in this paper
a property of measure dynamical systems, called
``admitting a separating sieve", which imitates the intrinsic
definition of topological distality.

\begin{defn}\label{defn-ssieve}
Let $\Xb$ be an ergodic dynamical system.
A sequence $A_1\supset A_2\supset \cdots$ of sets in $\Xcal$
with $\mu(A_n)>0$ and $\mu(A_n)\to 0$,
is called a {\bf separating sieve\/}
if there exists a
subset $X_0\subset X$ with $\mu(X_0)=1$ such that for every
$x,x'\in X_0$ the condition
``for every $n\in \N$ there exists $k\in \Z$ with
$T^k x,T^k x'\in A_n$" implies $x=x'$, or
in symbols:
$$
\bigcap_{n=1}^\infty\left(\bigcup_{k\in \Z}
T^k(A_n\times A_n)\right)\cap
(X_0 \times X_0)\subset \Del.
$$
We say that the ergodic system $\Xb$ is
{\bf measure distal\/} if either $\Xb$ is finite
or there exists a separating sieve.
\end{defn}

\br

Parry showed that every
measure dynamical system admitting a separating sieve
has zero entropy and that any
$T$-invariant measure on a minimal topologically distal
system gives rise to a measure dynamical system
admitting a separating sieve.

If $\Xb=(X,\Xcal,\mu,T)$ is an ergodic dynamical system and
$K\subset \Aut(\Xb)$ is a compact subgroup (where
$\Aut(\Xb)$ is endowed with the weak topology)
then the system $\Yb=\Xb/K$ is well defined and we say
that the extension $\pi:\Xb \to\Yb$ is a {\bf group extension\/}.
Using \eqref{iso-gr-diag} we can define the notion of
isometric extension or homogeneous extension in the
measure category. We will say that an ergodic system
{\bf admits a Furstenberg tower\/} if it is obtained
as a (necessarily countable) transfinite tower of
measure isometric extensions.
In 1976 in two outstanding papers \cite{Z1}, \cite{Z2}
R. Zimmer developed the theory of distal systems
(for a general locally compact acting
group). He showed that, as in the topologically distal case,
systems admitting Parry's separating sieve are exactly
those systems which admit Furstenberg towers.

\begin{thm}
\label{zimmer-structure-tm}
An ergodic dynamical system is measure distal iff
it admits a Furstenberg tower.
\end{thm}

\br

In \cite{Lis2} E. Lindenstrauss shows that
every ergodic measure distal $\Z$-system
can be represented as a minimal topologically distal
system. For the exact result see Theorem \ref{distal-model}
below.

\br

\section{Furstenberg-Zimmer structure theorem vs. its topological PI
version}\label{Sec-FZ}

Zimmer's theorem for distal systems leads directly to a
structure theorem for the general ergodic system.
Independently, and at about the same time,
Furstenberg proved the same theorem,
\cite{Fur4}, \cite{Fur5}.
He used it as the main tool for his proof of Szemer\'edi's theorem
on arithmetical progressions.
Recall that an extension
$\pi:(X,\Xcal,\mu,T)\to (Y,\Ycal,\nu,T)$ is a {\bf  weakly
mixing extension\/} if the relative product system
$\Xb\underset{\Yb}{\times}\Xb$ is ergodic.
(The system $\Xb\underset{\Yb}{\times}\Xb$ is defined by
the $T\times T$ invariant measure
$$
\mu\underset{\nu}{\times}\mu=\int_Y\mu_y\times\mu_y\,d\nu(y),
$$
on $X\times X$, where $\mu=\int_Y\mu_y\,d\nu(y)$
is the disintegration of $\mu$ over $\nu$.)

\begin{thm}[The Furstenberg-Zimmer structure theorem]\label{FZsth}
Let $\Xb$ be an ergodic dynamical system.
\begin{enumerate}
\item
There exists a maximal distal factor
$\phi:\Xb\to\Zb$ with
$\phi$ is a weakly mixing extension.
\item
This factorization is unique.
\end{enumerate}
\end{thm}

\br

{\begin{center}{$\divideontimes$}\end{center}}

\br

Is there a general structure theorem for minimal topological
systems?
Here, for the first time, we see a strong divergence
between the measure and the topological theories.
The culpability for this divergence is to be found in the notions of
proximality and proximal extension, which arise
naturally in the topological theory but do not appear
at all in the measure theoretical context. In building
towers for minimal systems we have to use two building
blocks of extremely different nature (isometric and proximal)
rather than one (isometric) in the measure category.
A pair of points $(x,x')\in X\times X$ is called {\bf proximal\/}
if it is not distal, i.e. if $\inf_{n\in \Z} d(T^nx,T^nx')=0$.
An extension $\pi:(X,T)\to (Y,T)$ is called {\bf proximal\/}
if every pair in $R_\pi$ is proximal. The next theorem was
developed gradually by several authors (Veech, Glasner-Ellis-Shapiro,
and McMahon, \cite{V}, \cite{EGS}, \cite{McM1}, \cite{V1}).
We need first to introduce some definitions.
We say that a minimal dynamical system \XT is {strictly \bf PI\/}
(proximal isometric) if it admits a tower consisting of
proximal and isometric extensions. It is called a {\bf PI
system\/} if there is a strictly PI minimal system
$(\tilde{X},T)$ and a proximal extension $\tet:\tilde{X}\to X$.
An extension $\pi:X \to Y$ is a {\bf RIC extension\/}
(relatively incontractible) if for every $n\in \N$
and every $y\in Y$ the set of almost periodic points in
$X_y^n=\pi^{-1}(y)\times\pi^{-1}(y)\times\dots\times
\pi^{-1}(y)$ ($n$ times) is dense. (A point is called
{\bf almost periodic\/} if its orbit closure is minimal.)
It can be shown that
a every isometric (and more generally, distal) extension
is RIC. Also every RIC extension is open.
Finally a homomorphism $\pi:X\to Y$ is called
{\bf topologically weakly mixing\/} if the dynamical system
$(R_\pi,T\times T)$ is topologically transitive.

\br

The philosophy in the next theorem is to regard proximal
extensions as `negligible' and then the claim is, roughly (i.e.
up to proximal extensions), that every minimal system is
a weakly mixing extension of its maximal PI factor.

\begin{thm}[Structure theorem for minimal systems]\label{minsth}
Given a metric minimal system $(X,T)$, there exists a
countable ordinal $\eta$ and a canonically defined
commutative diagram (the canonical PI-Tower)
\begin{equation}
\xymatrix
        {X \ar[d]_{\pi}             &
     X_0 \ar[l]_{\tilde{\theta_0}}
         \ar[d]_{\pi_0}
         \ar[dr]^{\sigma_1}         & &
     X_1 \ar[ll]_{\tilde{\theta_1}}
         \ar[d]_{\pi_1}
         \ar@{}[r]|{\cdots}         &
     X_{\nu}
         \ar[d]_{\pi_{\nu}}
         \ar[dr]^{\sigma_{\nu+1}}       & &
     X_{\nu+1}
         \ar[d]_{\pi_{\nu+1}}
         \ar[ll]_{\tilde{\theta_{\nu+1}}}
         \ar@{}[r]|{\cdots}         &
     X_{\eta}=X_{\infty}
         \ar[d]_{\pi_{\infty}}          \\
        pt                  &
     Y_0 \ar[l]^{\theta_0}          &
     Z_1 \ar[l]^{\rho_1}            &
     Y_1 \ar[l]^{\theta_1}
         \ar@{}[r]|{\cdots}         &
     Y_{\nu}                &
     Z_{\nu+1}
         \ar[l]^{\rho_{\nu+1}}          &
     Y_{\nu+1}
         \ar[l]^{\theta_{\nu+1}}
         \ar@{}[r]|{\cdots}         &
     Y_{\eta}=Y_{\infty}
    }
\nonumber
\end{equation}
where for each $\nu\le\eta, \pi_{\nu}$
is RIC, $\rho_{\nu}$ is isometric, $\theta_{\nu},
{\tilde\theta}_{\nu}$ are proximal extensions and
$\pi_{\infty}$ is RIC and topologically weakly mixing extension.
For a limit ordinal
$\nu ,\  X_{\nu}, Y_{\nu}, \pi_{\nu}$
etc. are the inverse limits (or joins) of
$ X_{\iota}, Y_{\iota}, \pi_{\iota}$ etc. for $\iota
< \nu$.
Thus $X_\infty$ is a proximal extension of $X$ and a RIC
topologically weakly mixing extension of the strictly PI-system
$Y_\infty$.
The homomorphism $\pi_\infty$ is an isomorphism (so that
$X_\infty=Y_\infty$) iff $X$ is a PI-system.
\end{thm}

\br

We refer to \cite{Gl3} for a review on structure theory
in topological dynamics.

\br

\section{Entropy: measure and topological}\label{Sec-ent}

\subsection{The classical variational principle}
For the definitions and the classical results concerning
entropy theory we refer to \cite{HKat}, Section 3.7 for
measure theory entropy and Section 4.4 for metric and
topological entropy.
The variational principle asserts that for a topological
$\Z$-dynamical system $(X,T)$ the topological entropy equals the
supremum of the measure entropies computed over all
the invariant probability measures on $X$. It was already
conjectured in the original paper of Adler, Konheim and McAndrew
\cite{AKM} where topological entropy was introduced;\
and then, after many stages
(mainly by Goodwyn, Bowen and Dinaburg; see for example \cite{DGS})
matured into a theorem in Goodman's paper \cite{Goodm}.

\begin{thm}[The variational principle]\label{var-princ}
Let \XT be a topological dynamical system, then
$$
\htop(X,T)={\sup}\{h_\mu:\mu\in M_T(X)\}
={\sup}\{h_\mu:\mu\in M^{\erg}_T(X)\}.
$$
\end{thm}

This classical theorem has had a tremendous
influence on the theory of dynamical systems
and a vast amount of literature ensued,
which we will not try to trace here (see \cite[Theorem 4.4.4]{HKat}).
Instead we would like to present a more
recent development.

\br

\subsection{Entropy pairs and UPE systems}

As we have noted in the introduction,
the theories of measurable dynamics (ergodic theory) and
topological dynamics exhibit a remarkable parallelism.
Usually one translates `ergodicity' as
`topological transitivity',`weak mixing' as `topological
weak mixing', `mixing' as `topological
mixing' and `measure distal' as `topologically distal'.
One often obtains this way parallel
theorems in both theories, though the methods of proof
may be very different.

What is then the topological analogue of being a K-system?
In \cite{Bla1} and  \cite{Bla2} F. Blanchard introduced a notion of
`topological $K$' for $\Z$-systems which he called UPE
(uniformly positive entropy).
This is defined as follows: a topological dynamical system $(X,T)$
is called a UPE system  if every open cover
of $X$ by two non-dense open sets $U$ and $V$ has
positive topological entropy. A local version of this
definition led to the concept of an entropy pair.
A pair $(x,x') \in X \times X,\ x\not=x'$ is an
entropy pair if for every open cover
$\Ucal=\{U,V\}$ of $X$, with
$x \in {\inte} (U^c)$ and $x' \in {\inte} (V^c)$,
the topological entropy $h(\Ucal)$
is positive.
The set of entropy pairs is denoted
by $E_X=E_{(X,T)}$ and it follows that the system $(X,T)$ is UPE
iff  $E_X= (X\times X)\setminus \Del$.
In general  $E^*=E_X\cup \Delta$ is
a $T\times T$-invariant closed symmetric and reflexive relation.
Is it also transitive? When the answer to this
question is affirmative then the quotient system
$X/E_X^*$ is the topological analogue of the
Pinsker factor. Unfortunately this need not
always be true even when $(X,T)$ is a minimal
system (see \cite{GW55} for a counter example).

The following theorem was proved in
Glasner and Weiss \cite{GW35}.

\begin{thm}\label{K>UPE}
If the compact system $(X,T)$ supports an
invariant measure $\mu$ for which the corresponding measure
theoretical system $(X,\Xcal,\mu,T)$ is a $K$-system,
then $(X,T)$ is UPE.
\end{thm}

Applying this theorem together with the Jewett-Krieger theorem
it is now possible to obtain a great variety
of strictly ergodic UPE systems.

Given a $T$-invariant probability measure $\mu$
on $X$, a pair $(x,x') \in X \times X,\ x\not=x'$
is called a $\mu$-entropy pair if for every Borel partition
$\al =\{Q,Q^c\}$ of $X$ with
$x \in {\inte}(Q)$ and $x' \in {\inte}(Q^c)$
the measure entropy $h_\mu(\al)$ is positive. This definition was
introduced by Blanchard, Host, Maass, Mart\'{\i}nez and Rudolph in
\cite{B-R} as a local generalization of Theorem \ref{K>UPE}.
It was shown in \cite{B-R} that
for every invariant probability measure $\mu$
the set $E_\mu$ of $\mu$-entropy pairs is contained
in $E_X$.

\begin{thm}\label{mu-subset-e}
Every measure entropy pair is a topological entropy pair.
\end{thm}

As in \cite{GW35} the main issue here is to understand
the, sometimes intricate, relation between the combinatorial
entropy $h_c(\Ucal)$ of a cover $\Ucal$ and the measure theoretical
entropy $h_\mu(\ga)$ of a measurable partition $\ga$ subordinate
to $\Ucal$.

\begin{prop}\label{pro-mu-subset-e}
Let $\Xb=(X,\Xcal,\mu,T)$ be a measure dynamical system.
Suppose $\Ucal=\{U,V\}$ is a measurable cover such that
every measurable two-set partition $\ga=\{H,H^c\}$ which
(as a cover) is finer than $\Ucal$ satisfies
$h_\mu(\ga)>0$; then $h_c(\Ucal)>0$.
\end{prop}

Since for a $K$-measure $\mu$ clearly every
pair of distinct points is in $E_\mu$, Theorem \ref{K>UPE}
follows from Theorem \ref{mu-subset-e}.
It was shown in \cite{B-R} that
when $(X,T)$ is uniquely ergodic the converse of Theorem
\ref{mu-subset-e} is also true:
$E_X=E_\mu$ for the unique invariant measure $\mu$ on $X$.

\br

\subsection{A measure attaining the topological entropy
of an open cover}

In order to gain a better understanding of the
relationship between measure entropy pairs and topological
entropy pairs one direction of a
variational principle for open covers (Theorem
\ref{mcvp} below) was proved in Blanchard, Glasner and Host \cite{BGH}.
Two applications of this
principle were given in \cite{BGH};\
(i) the construction, for a general system $(X,T)$,
of a measure $\mu\in M_T(X)$ with $E_X=E_\mu$,
and (ii) the proof that under a homomorphism
$\pi:(X,\mu,T)\to (Y,\nu,T)$ every entropy pair in
$E_\nu$ is the image of an entropy pair in $E_\mu$.

We now proceed with the statement and proof of this
theorem which is of independent interest.
The other direction of this variational principle
will be proved in the following subsection.

\begin{thm}\label{mcvp}
Let $(X,T)$ be a topological dynamical system,
and $\Ucal$ an open cover of $X$, then there exists a measure
$\mu\in M_T(X)$ such that $h_\mu(\al)\ge \htop(\Ucal)$ for
all Borel partitions $\al$ finer than $\Ucal$.
\end{thm}

\br

A crucial element of the proof of the
variational principle is a combinatorial lemma which we present
next.
We let
$\phi: [0,1]\to\R$ denote the function
$$
\phi(x)=-t\log t\quad {\text{ for\ }} 0<t\le 1 ;\  \phi(0)=0 \ .
$$
Let $\Lf=\{1,2,\dots,\ell\}$ be a
finite set, called the {\bf alphabet\/}; sequences
$\om=\om_1\ldots \om_n\in \Lf^n$, for $n\ge 1$, are called
{\bf words of length $n$ on the alphabet $\Lf$\/} .
Let $n$ and $k$ be two integers with $1\leq k\le n$.

For every word $\om$ of length $n$ and every word $\tet$ of length $k$
on the same alphabet, we denote by $p(\tet|\om)$ the frequency of
appearances of $\tet$ in $\om$,\ i.e.
$$
p(\tet|\om)=\frac{1}{n-k+1}
{\card}\big\{i :\  1\leq i\le n-k+1,\;
\om_i\om_{i+1}\ldots\om_{i+k-1}=
\tet_1\tet_2\ldots\tet_{k}\big\} \ .
$$
For every word $\om$ of
length $n$ on the alphabet $\Lf$, we let
$$
H_k(\om)=\sum_{\tet\in
\Lf^k}\phi\big(p(\tet|\om)\big) \ .
$$

\begin{lem} For every $h>0$, $\ep>0$, every integer
$k\ge 1$ and every sufficiently large integer $n$,
$$
{\card}\big\{
\om\in \Lf^n : H_k(\om) \le kh \big\} \le \exp\big(
n(h+\ep)\big) \ .
$$
\end{lem}\label{combin}
{\bf Remark. }It is equally true that, if $h \le \log ({\card} \Lf)$,
for sufficiently large $n$,
$$
{\card}\big\{ \om\in \Lf^n : H_k(\om) \leq kh \big\} \geq \exp\big(
n(h-\ep)\big) \ .
$$
We do not prove this inequality here,
since we have no use for it in the sequel.

\begin{proof}
{\bf The case $ k=1$.}

We have
\begin{equation}\label{H1}
{\card}\big\{ \om\in \Lf^n :\  H_1(\om) \le h \big\}
=\sum_{q \in K} \frac{n!}{q_1!\ldots q_{\ell}!}
\end{equation}
where $K$ is the set of $q=(q_1,\ldots ,q_{\ell})\in \N^{\ell}$
such that
$$
\sum_{i=1}^{\ell} q_i=n\ {\text{ and }}\ \sum_{i=1}^{\ell}
 \phi(\frac{q_i}{n})\le h \ .
$$
By Stirling's formula, there exist two
universal constants $c$ and
$c'$ such that
$$
c\big(\frac{m}{e})^m\sqrt m \le m!\,
\le c'\big({\frac{m}{e}})^m\sqrt m
$$
for every $m>0$.
From this we deduce the existence of a constant $C(\ell)$ such
that for every $ q\in K$,
$$
{\frac{n!}{q_1!\ldots q_{\ell}!}} \leq C(\ell)
\exp\big( n\sum_{i=1}^{\ell}
\phi(\frac{q_i}{n})\big)
\leq C({\ell})\exp(nh)\ .
$$
Now the sum \eqref{H1} contains at most $(n+1)^{\ell}$ terms;
so that we have
$$
{\card}\big\{ \om\in \Lf^n :\  H_1(\om) \le h \big\} \le
(n+1)^{\ell}C(\ell)\exp(nh)\le \exp\big( n(h+\ep)\big)
$$
for all sufficiently large $n$, as was to be proved.

{\bf The case $k>1$.}

For every word $\om$ of length $n\ge 2k$ on the alphabet
$\Lf$, and for
$0\le j<k$, we let $n_j$ be the integral part of $\frac{n-j}{k}$,
and $\om^{(j)}$ the word
$$
(\om_{j+1}\ldots\om_{j+k})\;
(\om_{j+k+1}\ldots\om_{j+2k})\;
\ldots\;
(\om_{j+(n_j-1)k+1}\ldots\om_{j+n_jk})
$$
of length $n_j$ on the
alphabet $B=\Lf^k$.

Let now $\tet$ be a word of length $k$ on the alphabet
$\Lf$; we also
consider $\tet$ as an element of $B$.
One easily verifies that, for every word $\om$ of length $n$ on the
alphabet $\Lf$,
$$
\big| p(\tet|\om)-\frac{1}{k}\sum_{j=0}^{k-1}
p(\tet|\om^{(j)})\big| \le \frac{k}{n-2k+1} \ .
$$
The function $\phi$ being uniformly continuous,
we see that for sufficiently large $n$,
and for every word $\om$ of length $n$ on $\Lf$,
$$
\sum_{\tet\in B} \left| \phi\big( p(\tet|\om)\big) -\phi\left(
\frac{1}{k}\sum_{j=0}^{k-1}p(\tet|\om^{(j)})\right) \right| <
\frac{\ep}{2}
$$
and by convexity of $\phi$,
$$
{\frac{1}{k}}\sum_{j=0}^{k-1}H_1(\om^{(j)})= \frac{1}{
k}\sum_{j=0}^{k-1}\sum_{\tet\in B} \phi\big(p(\tet|\om^{(j)})\big)
\leq \frac{\ep}{2} +
\sum_{\tet\in \Lf^k}\phi\big(p(\tet|\om)\big) =\frac{\ep}{2} +
H_k(\om).
$$
Thus, if $H_k(\om)\le kh$, there exists a $j$ such that
$H_1(\om^{(j)})\le \frac{\ep}{2} + kh$.

Now, given $j$ and a word $u$ of length
$n_j$ on the alphabet $B$, there
exist $\ell^{n-n_jk}\leq \ell^{2k-2}$ words $\om$ of length
$n$ on $\Lf$ such that $\om^{(j)}=u$. Thus for sufficiently large $n$,
by the first part of the proof,
\begin{align*}
{\card}\big\{ \om\in \Lf^n :  H_k(\om) \le kh \big\} & \le
{\ell}^{2k-2}\sum_{j=0}^{k-1}{\card}\big\{ u\in B^{n_j}:H_1(u) \le
\frac{\ep}{2} + kh\big\}\\
& \le {\ell}^{2k-2}\sum_{j=0}^{k-1}\exp\big(n_j(\ep+kh)\big)\\
& \le
{\ell}^{2k-2}k\exp\big(n(\frac{\ep}{k}+h)\big) \le
\exp\big(n(h+\ep\big))\ .
\end{align*}
\end{proof}

\br

Let $(X,T)$ be a compact dynamical system.
As usual we denote by $M_T(X)$ the set of
$T$-invariant probability measures on $X$, and by
$M_T^{\erg}(X)$ the subset of ergodic measures.

We say that a partition $\al$ is finer than a cover
$\Ucal$ when every atom of $\al$ is contained in an element of
$\Ucal$.
If $\al=\{A_1,\ldots,A_{\ell}\}$ is a partition of $X$,
$x\in X$ and $N\in\N$, we write $\om(\al,N,x)$ for the word of
length $N$ on the alphabet
$\Lf=\{1,\ldots,{\ell}\}$ defined by
$$
\om(\al,N,x)_n=i\quad {\text{if}}\quad T^{n-1}x\in A_i,\qquad
1\le n\le N\ .
$$

\br

\begin{lem}Let $\Ucal$ be a cover of $X$, $h=\htop(\Ucal)$,
$K\ge 1$ an integer, and $\{\al_l:1\le l\le K\}$ a finite
sequence of partitions of $X$, all finer than $\Ucal$.
For every $\ep>0$ and sufficiently large $N$,
there exists an $x\in X$ such that
$$
H_k\big(\om(\al_l,N,x)\big)
\ge k(h-\ep)
\ {\text  {for every}  }\ k,l\ {\text { with}  }\  1\leq k,\;l\le K.
$$
\end{lem}
\begin{proof}
One can assume that all the partitions
$\al_l$ have the same number of elements $\ell$ and we
let $\Lf=\{1,\ldots,\ell\}$.
For $1\le k \le K$ and $N\ge K$, denote
$$
\Om(N,k)=\{\om \in \Lf^N :\  H_k(\om) < k(h-\ep)\} \ .
$$
By Lemma \ref{combin}, for sufficiently large $N$
$$
{\card}(\Om(N,k))\leq \exp(N(h-\ep/2))\ {\text{ for all}}\ k\le K.
$$
Let us choose such an $N$ which moreover satisfies
$K^2 <\exp(N\ep/2)$.
For $1\le k,l\le K$, let
$$
Z(k,l)= \{x\in X :\om(\al_l, N,x)\in\Om(N,k)\} \ .
$$

The set $Z(k,l)$ is the union of ${\card}(\Om(N,k))$ elements of
$(\al_l)_0^{N-1}$. Now this partition is finer than the cover
$\Ucal_0^{N-1}$, hence $Z(k,l)$ is covered by
$$
{\card}(\Om(N,k))\leq\exp(N(h-\ep/2))
$$
elements of $\Ucal_0^{N-1}$.
Finally,
$$
\bigcup_{1\leq k,l\leq K} Z(k,l)
$$
is covered by
$K^2\exp(N(h-\ep/2))<\exp(Nh)$ elements of $\Ucal_0^{N-1}$.
As every subcover of $\Ucal_0^{N-1}$ has at least
$\exp(Nh)$ elements,
$$
\bigcup_{1\le k,l\le K} Z(k,l)\neq X.
$$
This completes the proof of the lemma.
\end{proof}  \br

\begin{proof}[Proof of theorem \ref{mcvp}]
Let $\Ucal=\{U_1,\ldots,U_\ell\}$ be an open cover of $X$.
It is clearly
sufficient to consider Borel partitions $\al$ of $X$ of the form
\begin{equation}\label{one*}
\al=\{ A_1,\ldots,A_\ell\}\ {\text{ with}}\ A_i\subset U_i
\  {\text{ for every} }\  i.
\end{equation}

\br

{\bf Step 1:\ }Assume first that $X$
is $0$-dimensional.

The family of partitions finer than
$\Ucal$, consisting of clopen sets and satisfying
\eqref{one*} is countable;
let $\{\al_l:l\ge 1\}$ be an enumeration of this family.
According to the previous lemma, there exists a sequence of integers
$N_K$ tending to
$+\infty$ and a sequence $x_K$ of elements of $X$ such that:
\begin{equation}\label{doub*}
H_k\big(\om(\al_l, N_K,x_K)\big) \ge k(h-\frac{1}{K}) \
{\text{ for
every} }\ 1\le k,\;l \le K.
\end{equation}
Write
$$
\mu_K=\frac{1}{N_K}\sum_{i=0}^{N_K-1} \del_{T^ix_K} \ .
$$
Replacing the
sequence $\mu_K$ by a subsequence
(this means replacing the sequence
$N_K$ by a subsequence, and the sequence $x_K$ by the
corresponding subsequence preserving the property \eqref{doub*}),
one can assume that the sequence
of measures $\mu_K$ converges weak$^*$ to a probability
measure $\mu$.
This measure $\mu$ is clearly $T$-invariant. Fix $k,l\ge 1$,
and let $F$ be an atom of the partition $(\al_l)_0^{k-1}$,
with name
$\tet\in\{1,\ldots, \ell\}^k$. For every $K$ one has
$$
\big|\mu_K(F)-p\big(\tet| \om(\al_l , N_K,x_K)\big)
\big| \le \frac{2k}{N_K}.
$$
Now as $F$ is clopen,
\begin{align*}
\mu(F)&
=\lim_{K\to\infty}\mu_K(F)=
\lim_{K\to\infty} p\big(\tet| \om(\al_l, N_K,x_K)\big)\;
\ {\text{hence}}\\
\phi(\mu(F) )
&=\lim_{K\to\infty} \phi\big(p\big(\tet| \om(\al_l, N_K,x_K)\big)
\big)
\end{align*}
and, summing over $\tet\in \{1,\ldots,\ell\}^k$, one gets
$$
H_\mu\big( (\al_l)_0^{k-1}\big)
=\lim_{K\to\infty} H_k\big( \om(\al_l , N_K,x_K)\big) \ge kh.
$$
Finally, by sending $k$ to infinity one obtains $h_\mu(\al_l)\ge h$.

Now, as $X$ is $0$-dimensional, the family of partitions
$\{\al_l\}$ is dense in the
collection of Borel partitions of $X$ satisfying \eqref{one*},
with respect to the
distance associated with $L^1(\mu)$. Thus, $h_\mu(\al)\ge h$
for every partition of this kind.

\br

{\bf Step 2:\ }The general case.

Let us recall a well known fact: there
exists a topological system $(Y,T)$, where $Y$ is $0$-dimensional,
and a continuous surjective map $\pi: Y\to X$ with
$\pi\circ T=T\circ \pi$.

(Proof : as $X$ is a compact metric space, it is easy to construct
a Cantor set
$K$ and a continuous surjective $f: K\to X$.
Put
$$
Y=\{ y \in K^{\Z}:\  f(y_{n+1})=T f(y_n)\
{\text{ for every }}\ n\in\Z\}
$$
and let $\pi: Y\to X$ be defined by $\pi(y)=f(y_0)$.

$Y$ is a closed subset of $K{^\Z}$ --- where the latter is
equipped with the product topology --- and is invariant under the
shift $T$ on $K^{\Z}$. It is
easy to check that $\pi$ satisfies the required conditions.)

Let $\Vcal=\pi^{-1}(\Ucal)=\{\pi^{-1}(U_1),
\ldots,\pi^{-1}(U_d)\}$
be the preimage of $\Ucal$ under $\pi$ ;
one has $\htop(\Vcal)=\htop(\Ucal)=h$. By the
above remark, there exists $\nu\in M(Y,T)$ such that
$h_\nu(\Qcal)\ge h$
for every Borel partition
$\Qcal$ of $Y$ finer than $\Vcal$. Let
$\mu=\nu\circ\pi^{-1}$ the measure which is the image of
$\nu$ under $\pi$.
One has $\mu\in M_T(X)$ and, for every Borel partition $\al$ of
$X$ finer than $\Ucal$, $\pi^{-1}(\al)$ is a Borel partition of
$Y$ which is finer than $\Vcal$ with
$$
h_\mu(\al)=h_\nu\big( \pi^{-1}(\al)\big)\ge h.
$$
This completes the proof of the theorem.
\end{proof}

\br

\begin{cor}\label{corol1} Let $(X,T)$ be a topological system,
$\Ucal$ an open cover of $X$ and $\al$
a Borel partition finer than $\Ucal$, then, there exists a
$T$-invariant ergodic measure $\mu$ on $X$ such that
$h_\mu(\al)\ge \htop(\Ucal)$.
\end{cor}
\begin{proof} By Theorem \ref{mcvp} there exists
$\mu\in M_T(X)$ with $h_{\mu}(\al)\ge \htop(\Ucal)$;
let $\mu=\int_\om \mu_\om\,dm(\om )$
be its ergodic decomposition. The corollary
follows from the formula
$$
\int h_{\mu_\om}(\al)\,dm(\om) =h_\mu(\al).
$$
\end{proof}

\br

\subsection{The variational principle for open covers}
Given an open cover $\Ucal$ of the dynamical system \XT,
the results of the previous subsection imply the inequality
$$
\sup_{\mu\in M_T(X)}\inf_{\al\succ\Ucal} h_\mu(\al) \ge
\htop(\Ucal).
$$
We will now present a new result which will
provide the fact that
$$
\sup_{\mu\in M_T(X)}\inf_{\al\succ\Ucal} h_\mu(\al) =
\htop(\Ucal)
$$
thus completing the proof of a variational
principle for $\Ucal$.

\br

We first need a universal version of the Rohlin lemma.

\begin{prop}\label{univ-roh}
Let $(X,T)$ be a (Polish) dynamical system
and assume that there exists on $X$ a $T$-invariant
aperiodic probability measure. Given a positive
integer $n$ and a real number $\delta>0$ there exists
a Borel subset $B\subset X$ such that the sets
$B, TB, \dots, T^{n-1}B$ are pairwise disjoint and for every
aperiodic $T$-invariant probability measure $\mu\in M_T(X)$
we have $\mu(\bigcup_{j=0}^{n-1} T^j B) > 1 - \delta$.
\end{prop}

\begin{proof}
Fix $N$ (it should be larger than $n/\delta$ for the required height
$n$ and error $\delta$). The set of points that are periodic with
period $\le N$ is closed. Any point in the complement
(which by our assumption is nonempty) has,
by continuity, a neighborhood $U$ with $N$ disjoint forward iterates.
There is a countable subcover
$\{U_m\}$ of such sets since the space is Polish.
Take $A_1=U_1$ as a base for a {\em Kakutani sky-scraper}
\begin{gather*}
\{T^j A_1^k: j=0,\dots, k-1;\ k=1,2,\dots\},\\
A_1^k=\{x\in A_1: r_{A_1}(x)=k\},
\end{gather*}
where $r_{A_1}(x)$ is the first integer $j\ge 1$ with $T^jx\in A_1$.
Next set
$$
B_1= \bigcup_{k\ge 1}\bigcup_{j=0}^{[(k-n-1)/n]} T^{jn}A_1^k,
$$
so that the sets $B_1, TB_1,\dots, T^{n-1}B_1$ are pairwise disjoint.

Remove the full forward $T$ orbit of $U_1$ from the space and repeat
to find $B_2$ using as a base for the next Kakutani
sky-scraper $A_2$ defined as $U_2$ intersected with the part of $X$ not
removed earlier. Proceed by induction to define the
sequence $B_i, \ i=1,2,\dots$ and set $B=\bigcup_{i=1}^\infty B_i$.
By Poincar\'e recurrence for any aperiodic invariant measure we exhaust
the whole space except for $n$ iterates of the union $A$ of the bases
of the Kakutani sky-scrapers.
By construction $A=\bigcup_{m=1}^\infty A_m$ has $N$ disjoint iterates
so that $\mu(A) \le 1/N$ for every $\mu\in M_T(X)$.
Thus $B, TB,\dots, T^{n-1}B$ fill all but $n/N < \delta$
of the space uniformly
over the aperiodic measures $\mu\in M_T(X)$.
\end{proof}

\br

Let \XT be a dynamical system and $\Ucal=\{U_1,U_2,
\dots U_{\ell}\}$ a finite open cover.
We denote by $\Acal$ the collection of all finite Borel
partitions $\al$ which refine $\Ucal$, i.e. for every
$A\in \al$ there is some $U\in \Ucal$ with $A\subset U$.
We set
$$
\check{ h}(\Ucal)=
\sup_{\mu\in M_T(X)}\inf_{\al\in\Acal} h_\mu(\al) \qquad
{\text{and}} \qquad
\hat h(\Ucal)=
\inf_{\al\in\Acal}\sup_{\mu\in M_T(X)} h_\mu(\al).
$$

\begin{prop}\label{bar-hat}
Let $(X,T)$ be a dynamical system, $\Ucal=\{U_1,U_2,
\dots U_{\ell}\}$ a finite open cover, then
\begin{enumerate}
\item
$\check{ h}(\Ucal) \le \hat h(\Ucal)$,
\item
$\hat h(\Ucal) \le \htop(\Ucal)$.
\end{enumerate}
\end{prop}

\begin{proof}
1.\
Given $\nu\in M_T(X)$ and $\al\in \Acal$ we obviously have
$h_\nu(\al) \le \sup_{\mu\in M_T(X)} h_\mu(\al)$. Thus
$$
\inf_{\al\in\Acal}h_\nu(\al) \le \inf_{\al\in\Acal}
\sup_{\mu\in M_T(X)} h_\mu(\al) = \hat h(\Ucal),
$$
and therefore also $\check{ h}(\Ucal) \le \hat h(\Ucal)$.

2.\
Choose for $\ep> 0$ an integer $N$ large enough so that
there is a subcover $\Dcal \subset \Ucal_0^{N-1}=
\bigvee_{j=0}^{N-1} T^{-j}\Ucal$ of cardinality
$2^{N(\htop(\Ucal)+\ep)}$. Apply Proposition \ref{univ-roh}
to find a set $B$ such that  the sets
$B, TB, \dots ,T^{N-1}B$ are pairwise disjoint and for every
$T$-invariant Borel probability measure $\mu\in M_T(X)$
we have $\mu(\bigcup_{j=0}^{N-1} T^j B) > 1 - \delta$.
Consider $\Dcal_B=\{D\cap B: D\in \Dcal\}$,
the restriction of the cover $\Dcal$ to $B$, and find
a partition $\beta$ of $B$ which refines $\Dcal_B$.
Thus each element $P\in \beta$ has the form
$$
P=P_{i_0,i_1,\dots,i_{N-1}} \subset
\left(\bigcap_{j=0}^{N-1} T^{-j} U_{i_j}\right)
\cap B,
$$
where $\bigcap_{j=0}^{N-1} T^{-j} U_{i_j}$ represents a
typical element of $\Dcal$.
Next use the partition $\beta$ of $B$ to define a partition
$\al=\{A_i:i=1,\dots,\ell\}$ of $\bigcup_{j=0}^{N-1} T^{j}B$
by assigning to the set $A_i$ all sets of the form
$T^jP_{i_0,i_1,\dots,i_j,\dots,i_{N-1}}$ where $i_j=i$
($j$ can be any number in $[0,N-1]$).
On the remainder of the space $\al$ can be taken to be any
partition refining $\Ucal$.

Now if $N$ is large and $\del$ small
enough then
\begin{equation}\label{estim}
h_\mu(\al) \le \htop(\Ucal)+2\ep.
\end{equation}
Here is a sketch of how one establishes this inequality.
For $n >> N$ we will estimate $H_\mu(\al_0^{n-1})$ by
counting how many $(n,\al)$-names are needed to cover
most of the space. We take $\del>0$ so that $\sqrt{\del}<< \ep$.
Denote $E = B \cup TB \cup \dots \cup T^{N-1}B$
(so that $\mu(E) > 1 - \delta$).
Define
$$
f(x)=\frac{1}{n}\sum_{i=0}^n \ch_E(T^ix),
$$
and observe that $0 \le f \le1$
and
$$
\int_X f(x)\, d\mu(x) > 1 - \del,
$$
since $T$ is measure preserving.
Therefore $\int (1-f) < \del$ and (Markov's
inequality)
$$
\mu \{x : (1-f) \ge \sqrt{\del}\}
\le \frac{1}{\sqrt{\del}} \int(1-f) \le \sqrt{\del}.
$$
It follows that for points $x$ in $G =\{f > 1- \sqrt{\del}\}$,
we have the property that $T^i x \in E$ for most $i$ in $[0,n]$.

Partition $G$ according to the values of $i$ for which $T^i x\in B$.
This partition has at most
$$
\sum_{j\le\frac{n}{N}}\binom{n}{j} \le \frac{n}{N}\binom{n}{n/N}
$$
sets, a number which is exponentially small in $n$
(if $N$ is sufficiently large).

For a fixed choice of these values the times when we are
not in $E$ take only $n\sqrt{\del}$ values and there we have
$ < l^{n\sqrt{\del}}$ choices.

Finally when $T^i x \in B$ we have at most
$2^{(N({\htop}(U) + \epsilon))}$
names so that the total contribution is $<
2^{(N({\htop}(U) + \epsilon))\frac{n}{N}}$.

Collecting these estimations we find that
$$
H(\alpha_0^{n-1}) < n({\htop}(U) + 2\epsilon),
$$
whence \eqref{estim}.
This completes the proof of the proposition.
\end{proof}

\br

We finally obtain:


\begin{thm}[The variational principle for open covers]\label{cvp}
Let $(X,T)$ be a dynamical system, $\Ucal=\{U_1,U_2,
\dots U_k\}$ a finite open cover
and denote by $\Acal$ the collection of all finite Borel
partitions $\al$ which refine $\Ucal$, then
\begin{enumerate}
\item
for every $\mu \in M_T(X)$,
$\inf_{\al\in\Acal} h_\mu(\al) \le \htop(\Ucal)$,
and
\item
there exists an ergodic measure $\mu_0\in M_T(X)$
with  $h_{\mu_0}(\al) \ge  \htop(\Ucal)$ for every
Borel partition $\al \in \Acal$.
\item
$$
\check{ h}(\Ucal) = \hat h(\Ucal) =  \htop(\Ucal).
$$
\end{enumerate}
\end{thm}

\begin{proof}
1.\
This assertion can be formulated by the inequality
$\check{ h}(\Ucal) \le \htop(\Ucal)$ and it
follows by combining the two parts of Lemma \ref{bar-hat}.

2.\
This is the content of Theorem \ref{mcvp}.

3.\
Combine assertions 1 and 2.
\end{proof}

\br

\subsection{Further results connecting topological
and measure entropy}

Given a topological dynamical system \XT and a measure
$\mu\in M_T(X)$,
let $\pi:(X,\Xcal,\mu,T)\to (Z,\Zcal,\eta,T)$ be the
{\bf measure-theoretical\/} Pinsker factor of $(X,\Xcal,\mu,T)$, and let
$\mu=\int_Z \mu_z\,d\eta(z)$ be the disintegration of
$\mu$ over $(Z,\eta)$. Set
$$
\lambda=\int_Z (\mu_z \times \mu_z)\,d\eta(z),
$$
the relatively independent joining of $\mu$ with itself over $\eta$.
Finally let $\Lambda_\mu={\supp}(\lambda)$
be the topological support of $\lambda$ in $X\times X$.
Although the Pinsker factor is, in general, only defined
measure theoretically, the measure $\la$ is a well defined
element of $M_{T\times T}(X\times X)$.
It was shown in Glasner \cite{Gl25} that $E_\mu=
\La_\mu\setminus \Del$.

\begin{thm}\label{char-mu-ep}
Let \XT be a topological dynamical system and let
$\mu\in\linebreak[0] M_T(X)$.
\begin{enumerate}
\item
$E_\mu=\Lambda_\mu\setminus \Delta$ and
$\La_\mu=E_\mu\cup \{(x,x): x\in {\supp}(\mu)\}$.
\item
${\cls} E_\mu\subset \Lambda_\mu$.
\item
If $\mu$ is ergodic with positive entropy then
${\cls} E_\mu=\Lambda_\mu$.
\end{enumerate}
\end{thm}

\br

One consequence of this characterization of the set of
$\mu$-entropy pairs is a description of the set
of entropy pairs of a product system.
Recall that an $E$-system is a system for which there exists
a probability invariant measure with full support.

\begin{cor}\label{prod-ep}
Let $(X_1,T)$ and $(X_2,T)$ be two topological $E$-systems then:
\begin{enumerate}
\item
$E_{X_1\times X_2}=
(E_{X_1}\times E_{X_2}) \cup (E_{X_1}\times \Del_{X_2})
\cup (\Del_{X_1}\times E_{X_2})$.
\item
The product of two UPE systems is UPE.
\end{enumerate}
\end{cor}

\br

Another consequence is:

\begin{cor}\label{ep-proximal}
Let $(X,T)$ be a topological dynamical system,
$P$ the proximal relation on $X$. Then:
\begin{enumerate}
\item
For every $T$-invariant ergodic measure $\mu$
of positive entropy
the set $P\cap E_\mu$ is residual in the
$G_\delta$ set $E_\mu$ of $\mu$ entropy pairs.
\item
When $E_X\ne\emptyset$ the set $P\cap E_X$ is residual in the
$G_\delta$ set  $E_X$ of topological entropy pairs.
\end{enumerate}
\end{cor}

\br

Given a dynamical system \XT,
a pair $(x,x')\in X\times X$ is called
a {\bf Li--Yorke pair\/} if it is a proximal pair but not
an asymptotic pair.
A set $S\subseteq X$ is called {\bf scrambled\/} if any pair of
distinct points $\{x,y\}\subseteq S$ is a Li--Yorke pair. A dynamical
system $(X,T)$ is called {\bf chaotic in the sense
of Li and Yorke\/} if there is an uncountable scrambled set.
In \cite{BGKM} Theorem \ref{char-mu-ep} is applied to solve
the question whether positive topological entropy implies
Li--Yorke chaos as follows.

\begin{thm}
Let $(X,T)$ be a topological dynamical system.
\begin{enumerate}
\item
If $(X,T)$ admits
a $T$-invariant ergodic measure $\mu$ with respect to which the measure
preserving system $(X,\Xcal,\mu,T)$ is not measure distal then $(X,T)$
is Li--Yorke chaotic.
\item
If $(X,T)$ has positive topological entropy then it is
Li--Yorke chaotic.
\end{enumerate}
\end{thm}

\br

In \cite{BHR} Blanchard, Host and Ruette show that
in positive entropy systems there are also many
asymptotic pairs.

\begin{thm}\label{BHR}
Let $(X,T)$ be a topological dynamical system with positive
topological entropy. Then
\begin{enumerate}
\item
The set of points $x\in X$ for which there is
some $x'\ne x$ with $(x,x')$ an asymptotic pair, has
measure $1$ for every invariant probability measure on $X$
with positive entropy.
\item
There exists a probability measure $\nu$ on
$X\times X$ such that $\nu$ a.e. pair $(x,x')$
is Li--Yorke and positively asymptotic; or more precisely
for some $\delta >0$
\begin{gather*}
\lim_{n\to+\infty}d(T^nx,T^nx')=0, \qquad {\text and}\\
\liminf_{n\to+\infty}d(T^{-n}x,T^{-n}x')=0,
\qquad
\limsup_{n\to+\infty}d(T^{-n}x,T^{-n}x')\ge\delta.
\end{gather*}
\end{enumerate}
\end{thm}

\br

\subsection{Topological determinism and zero entropy}

Following \cite{KSS} call a dynamical system \XT
{\bf deterministic\/} if every $T$-factor is
also a $T^{-1}$-factor. In other words
every closed equivalence relation $R\subset
X\times X$ which has the property $TR\subset R$
also satisfies $T^{-1}R\subset R$.
It is not hard to see that an equivalent condition is
as follows. For every continuous real valued function
$f\in C(X)$ the function $f\circ T^{-1}$ is
contained in the smallest closed subalgebra $\Acal\subset C(X)$
which contains the constant function $\ch$ and the
collection $\{f\circ T^n: n\ge 0\}$.
The folklore question whether the latter condition implies zero
entropy was open for awhile. Here we note that
the affirmative answer is a direct consequence of
Theorem \ref{BHR} (see also \cite{KSS}).

\begin{prop}\label{non-inv-f}
Let $(X,T)$ be a topological dynamical system
such that there exists a $\delta>0$ and a pair
$(x,x')\in X\times X$ as in
Theorem \ref{BHR}.2. Then \XT is not deterministic.
\end{prop}

\begin{proof}
Set
$$
R=\{(T^nx,T^nx'): n\geq 0\}
\cup \{(T^nx',T^nx): n\geq 0\}
\cup \Delta.
$$
Clearly $R$ is a closed equivalence relation
which is $T$-invariant but not $T^{-1}$-invariant.
\end{proof}

\begin{cor}
A topologically deterministic dynamical system has zero entropy.
\end{cor}

\begin{proof}
Let $(X,T)$ be a topological dynamical system with positive
topological entropy; by Theorem \ref{BHR}.2.
and Proposition \ref{non-inv-f} it is not deterministic.
\end{proof}

\br

%

\newpage


{\large{\part{Meeting grounds}}}

\br

\section{Unique ergodicity}\label{Sec-JK}

The topological system
$(X,T)$ is called {\bf uniquely ergodic\/} if $M_T(X)$
consists of a single element $\mu$. If in addition $\mu$ is
a full measure (i.e. ${\supp} \mu=X$) then the system
is called {\bf strictly ergodic\/}
(see \cite[ Section 4.3]{HKat}).
Since the ergodic measures are characterized as the extreme points
of the Choquet simplex $M_T(X)$, it follows immediately
that a uniquely ergodic measure is ergodic.
For a while it was believed that strict ergodicity ---
which is known to imply some strong topological consequences
(like in the case of $\Z$-systems, the fact that {\em every\/}
point of $X$ is a generic point and moreover that the
convergence of the ergodic sums $\A_n(f)$
to the integral $\int f\, d\mu, \ f\in C(X)$
is {\em uniform\/}) ---
entails some severe restrictions on the measure-theoretical
behavior of the system. For example, it was believed that
unique ergodicity implies zero entropy. Then, at first some examples
were produced to show that this need not be the case.
Furstenberg in \cite{Fur3} and Hahn and Katznelson
in \cite{HK} gave examples of uniquely ergodic systems with
positive entropy. Later in 1970 R. I. Jewett surprised everyone
with his outstanding result: every weakly mixing measure
preserving $\Z$-system has a strictly ergodic model,
\cite{Jew}.
This was strengthened by Krieger \cite{K} who showed that
even the weak mixing assumption is redundant and that
the result holds for every ergodic $\Z$-system.

We recall the following well known characterizations
of unique ergodicity (see \cite[Theorem 4.9]{G}).

\begin{prop}\label{unique-erg}
Let \XT be a topological system. The following conditions are
equivalent.
\begin{enumerate}
\item
\XT is uniquely ergodic.
\item
$C(X)=\R+\bar B$, where $B=\{g-g\circ T: g\in C(X)\}$.
\item
For every continuous function $f\in C(X)$
the sequence of functions
$$
\A_nf(x)=\frac{1}{n}\sum_{j=0}^{n-1}f(T^jx).
$$
converges uniformly to a constant function.
\item
For every continuous function $f\in C(X)$
the sequence of functions $\A_n(f)$ converges pointwise to a
constant function.
\item
For every function $f\in A$, for a collection
$A\subset C(X)$ which linearly spans a uniformly dense
subspace of $C(X)$, the sequence of functions $\A_n(f)$
converges pointwise to a constant function.
\end{enumerate}
\end{prop}

\br

Given an ergodic dynamical system
$\Xb=(X,\Xcal,\mu,T)$
we say that the system $\hat \Xb=(\hat X,\hat \Xcal,\hat \mu,T)$ is
a {\bf topological model\/} (or just a model) for
$\Xb$ if $(\hat X,T)$ is a topological system,
$\hat \mu\in M_T(\hat X)$ and the systems
$\Xb$ and $\hat \Xb$ are measure theoretically isomorphic.
Similarly we say that $\hat \pi:\hat \Xb \to \hat \Yb$
is a {\bf topological model\/} for $\pi:\Xb\rightarrow \Yb$
when $\hat \pi$ is a topological factor map and there
exist measure theoretical isomorphisms $\phi$
and $\psi$ such that the diagram
\begin{equation*}\label{eq-model2}
\xymatrix
{
\Xb \ar[d]_{\pi} \ar[r]^{\phi}  &
  \hat \Xb \ar[d]^{\hat \pi} \\
\Yb\ar[r]_{\psi} &  \hat\Yb
}
\end{equation*}
is commutative.

\br

\section{The relative Jewett-Krieger theorem}
In this subsection we will prove the following generalization
of the Jewett-Krieger theorem (see \cite[Theorem 4.3.10]{HKat}).

\begin{thm}\label{rel-jk}
If $\pi:\Xb=(X,\Xcal,\mu,T)\rightarrow \Yb=(Y,\Ycal,\nu,T)$ is
a factor map with $\Xb$ ergodic and $\hat\Yb$ is a uniquely
ergodic model for $\Yb$ then there
is a uniquely ergodic model $\hat \Xb$ for $\Xb$ and a factor map
$\hat \pi:\hat \Xb \to \hat \Yb$  which is a model for
$\pi:\Xb\rightarrow \Yb$.
\end{thm}

In particular, taking $\Yb$ to be the trivial one point system we
get:

\begin{thm}\label{jk}
Every ergodic system has a uniquely ergodic model.
\end{thm}

\br

Several proofs have been given of this theorem, e.g.
see \cite{DGS} and \cite{BF}.
We will sketch a proof which will serve the relative case as well.

\begin{proof}[Proof of theorem \ref{rel-jk}]

A key notion for this proof is that of a {\bf uniform} partition
whose importance in this context was emphasized by
G. Hansel and J.-P. Raoult, \cite{HR}.
\begin{defn}
A set $B \in \Xcal$ is uniform if
$$
\lim_{N \to \infty}\ {\esssup}_x
\left|\ \frac 1N\ \sum^{N-1}_0\ 1_B(T^ix)- \mu(B)\right|=0.
$$
A partition $\Pcal$ is uniform if, for all $N$, every set in
$\bigvee^N_{-N}\ T^{-i}\Pcal$ is uniform.
\end{defn}

The connection between uniform sets, partitions and unique
ergodicity lies in Proposition \ref{unique-erg}.
It follows easily from that
proposition that if $\Pcal$ is a uniform partition, say into the
sets $\{P_1,\ P_2, \ldots, P_a\}$, and we denote by $\Pcal$ also
the mapping that assigns to $x \in X$, the index $1 \leq i \leq a$
such that $x \in P_i$, then we can map $X$ to $\{1,\ 2, \ldots,
a\}^\mathbb Z = A^\mathbb Z$ by:
$$
\pi(x)= (\ldots, \Pcal (T^{-1}x),\ \Pcal(x),\ \Pcal(Tx), \ldots,
\Pcal(T^nx), \ldots).
$$
Pushing forward the measure $\mu$ by $\pi$, gives $\pi \circ \mu$
and the closed support of this measure will be a closed shift
invariant subset, say $E \subset A^\mathbb Z$.  Now the indicator
functions of finite cylinder sets span the continuous functions on
$E$, and the fact that $\Pcal$ is a uniform partition and
Proposition \ref{unique-erg} combine to establish that
$(E,\ $shift) is uniquely
ergodic.  This will not be a model for $(X,\ \Xcal,\ \mu,\ T)$
unless $\bigvee^\infty_{- \infty}\ T^{-i}\Pcal= \Xcal$
modulo null sets, but in any case this does give a model for a
nontrivial factor of $X$.

Our strategy for proving Theorem \ref{jk} is to first construct a single
nontrivial uniform partition.  Then this partition will be refined
more and more via uniform partitions until we generate the entire
$\sigma$-algebra $\Xcal$.  Along the way we will be showing how
one can prove a relative version of the basic Jewett--Krieger
theorem.  Our main tool is the use of Rohlin towers.  These are
sets $B \in \Xcal$ such that for some $N,\ B,\ TB, \ldots,
T^{N-1}B$ are disjoint while $\bigcup^{N-1}_0\ T^iB$ fill up
most of the space. Actually we need Kakutani--Rohlin towers, which
are like Rohlin towers but fill up the whole space.  If the
transformation does not have rational numbers in its point spectrum
this is not possible with a single height, but two heights that are
relatively prime, like $N$ and $N+1$ are certainly possible.  Here
is one way of doing this. The ergodicity of $(X,\ \Xcal,\ \mu,\
T)$ with $\mu$ non atomic easily yields, for any $n$, the existence
of a positive measure set $B$, such that
$$
T^i\ B \cap B = \emptyset \qquad , \qquad i=1,\ 2, \ldots, n.
$$
With $N$ given, choose $n \geq 10 \cdot N^2$ and find $B$ that
satisfies the above.  It follows that the return time
$$
r_B(x)= \inf\{i>0:T^ix \in B\}
$$
is greater than $10 \cdot N^2$ on $B$.  Let
$$
B_\ell = \{x:r_B(x)=\ell\}.
$$
Since $\ell$ is large (if $B_\ell$ is nonempty) one can write
$\ell$ as a positive combination of $N$ and $N \times 1$, say
$$
\ell = Nu_\ell + (N+1) v_\ell.
$$
Now divide the column of sets $\{T^iB_\ell : 0 \leq i < \ell\}$
into $u_\ell$-blocks of size $N$ and $v_\ell$-blocks of size $N+1$
and mark the first layer of each of these blocks as belonging to
$C$.  Taking the union of these marked levels ($T^iB_\ell$ for
suitably chosen $i$) over the various columns gives us a set $C$
such that $r_C$ takes only two values -- either $N$ or $N+1$ as
required.

It will be important for us to have at our disposal K-R towers like
this such that the columns of say the second K-R tower are composed
of entire subcolumns of the earlier one.  More precisely we want
the base $C_2$ to be a subset of $C_1$ -- the base of the first
tower. Although we are not sure that this can be done with just two
column heights we can guarantee a bound on the number of columns
that depends only on the maximum height of the first tower.  Let us
define formally:
\begin{defn}
A set $C$ will be called the base of
a {\bf bounded} K-R tower if for some
$N,\ \bigcup^{N-1}_0\ T^iC=X$ up to a $\mu$-null set.
The least $N$ that satisfies this will be called the
{\bf height} of $C$,
and partitioning $C$ into sets of constancy of $r_C$ and viewing
the whole space $X$ as a tower over $C$ will be called the
K-R tower with columns the sets $\{T^iC_\ell:\ 0 \leq i <\ell\}$
for $C_\ell=\{x \in C: r_C(x)=\ell\}$.
\end{defn}

Our basic lemma for nesting these K-R towers is:
\begin{lem}\label{4N}
Given a bounded K-R tower with base $C$ and height $N$, for
any $n$ sufficiently large there is a bounded K-R tower with base
$D$ contained in $C$ whose column heights are all at least $n$ and
at most $n+4N$.
\end{lem}

\begin{proof}
We take an auxiliary set $B$ such that $T^i\ B \cap
B= \emptyset$ for all $0 <i<10(n+2N)^2$ and look at the unbounded
(in general) K-R tower over $B$.  Using the ergodicity it is easy
to arrange that $B \subset C$.  Now let us look at a single column
over $B_m$, with $m \geq 10\ (n+2N)^2$.  We try to put down blocks
of size $n+2N$ and $n+2N+1$, to fill up the tower. This can
certainly be done but we want our levels to belong to $C$. We can
refine  the column over $B_m$ into a finite number of columns so
that each level is either entirely within $C$ or in $X \backslash
C$.  This is done by partitioning the base $C$ according to the
finite partition:
$$
\bigcap^{m-1}_{i=0}\ T^{-1}\{C,\ X \backslash C\}.
$$
Then we move the edge of each block to the nearest level that
belongs to $C$.  The fact that the height of $C$ is $N$ means that
we do not have to move any level more than $N-1$ steps, and so at
most we lose $2N-2$ or gain that much thus our blocks, with bases
now all in $C$, have size in the interval $[n,\ n+4N]$ as
required.
\end{proof}

It is clear that this procedure can be iterated to give an infinite
sequence of nested K-R towers with a good control on the variation
in the heights of the columns.  These can be used to construct
uniform partitions in a pretty straightforward way, but we need one
more lemma which strengthens slightly the ergodic theorem.  We will
want to know that when we look at a bounded K-R tower with base $C$
and with  minimum column height sufficiently large that for most of
the fibers of the towers (that is for $x \in C,\ \{T^ix:\ 0 \leq
i<r_C(x)\}$) the ergodic averages of some finite set of functions
are close to the integrals of the functions.  It would seem that
there is a problem because the base of the tower is a set of very
small measure (less than 1/min column height) and it may be that
the ergodic theorem is not valid there.  However, a simple
averaging argument using an intermediate size gets around this
problem.  Here is the result which we formulate for simplicity for
a single function $f$:

\begin{lem}\label{erg}
Let $f$ be a bounded
function and $(X,\ \Xcal,\ \mu,\ T)$ ergodic.  Given $\epsilon
>0$, there is an $n_0$, such that if a bounded K-R tower with base
$C$ has minimum column height at least $n_0$, then those fibers
over $x \in C:\ \{T^ix:\ 0 \leq i < r_C(x)\}$ that satisfy
$$
\left| \frac {1}{r_C(x)}\ \sum^{r_C(x)-1}_{i=0}\
f(T^ix)-\int_X\ fd\mu \right| <\epsilon
$$
fill up at least $1-\epsilon$ of the space.
\end{lem}

\begin{proof}
Assume without loss of generality that $|f| \le 1$.  For
a $\delta$ to be specified later find an $N$ such that the set of
$y \in X$ which satisfy
\begin{equation}\label{*}
\left|\frac 1N\ \sum^{N-1}_0\ f(T^iy)-\int fd \mu\right| <
\delta
\end{equation}
has measure at least $1-\delta$.  Let us denote the set of $y$ that
satisfy \eqref{*} by $E$.  Suppose now that $n_0$ is large enough so that
$N/n_0$ is negligible -- say at most $\delta$.  Consider a bounded
K-R tower with base $C$ and with minimum column height greater than
$n_0$. For each fiber of this tower, let us ask what is the
fraction of its points that lie in $E$.  Those fibers with at least
a $\sqrt{\delta}$ fraction of its points not in $E$ cannot fill up
more than a $\sqrt{\delta}$ fraction of the space, because
$\mu(E)>1-\delta$.

Fibers with more than $1-\sqrt{\delta}$ of its points
lying in $E$ can be divided into disjoint blocks of size $N$
that cover all the points that lie in $E$.
This is done by starting at $x\in C$, and moving up the fiber,
marking the first point in $E$, skipping $N$ steps and continuing
to the next point in $E$ until we exhaust the fiber.
On each of these $N$-blocks the average of $f$ is within $\del$
of its integral, and since $|f|\le 1$ if $\sqrt{ \del } <
\ep/10$ this will guarantee that the average of $f$ over the
whole fiber is within $\ep$ of its integra.
\end{proof}

We are now prepared to construct uniform partitions.  Start with
some fixed nontrivial partition $\Pcal_0$.  By Lemma \ref{erg},
for any tall enough bounded K-R tower at least 9/10 of the columns will
have the 1-block distribution of each $\Pcal_0$-name within $\frac
{1}{10}$ of the actual distribution.  We build a bounded K-R tower
with base $C_1(1)$ and heights $N_1,\ N_1+1$ with $N_1$ large
enough for this to be valid.  It is clear that we can modify
$\Pcal_0$ to $\Pcal_1$ on the bad fibers so that now all fibers have a
distribution of 1-blocks within $\frac {1}{10}$ of a fixed
distribution.  We call this new partition $\Pcal_1$.  Our
further changes in $\Pcal_1$ will not change the $N_1,\ N_1+1$
blocks that we see on fibers of a tower over our ultimate $C_1$.
Therefore, we will get a uniformity on all blocks of size $100N_1$.
The 100 is to get rid of the edge effects since we only know the
distribution across fibers over points in $C_1(1)$.

Next we apply Lemma \ref{erg} to the 2-blocks in $\Pcal_1$ with 1/100.
We choose $N_2$ so large that $N_1/N_2$ is negligible and so that
any uniform K-R tower with height at least $N_2$ has for at least
99/100 of its fibers a distribution of 2-blocks within $1 / 100$ of
the global $\Pcal_1$ distribution.
Apply Lemma \ref{4N} to find a uniform K-R tower
with base $C_2(2) \subset C_1(1)$ such that its column heights are
between $N_2$ and $N_2+4N_1$.  For the fibers with good $\Pcal_1$
distribution we make no change.  For the others, we copy on most of
the fiber (except for the top $10 \cdot N^2_1$ levels) the
corresponding $\Pcal_1$-name from one of the good columns.  In
this copying we also copy the $C_1(1)$-name so that we preserve the
blocks. The final $10 \cdot N^2_1$ spaces are filled in with $N_1,\
N_1+1$ blocks.  This gives us a new base for the first tower that
we call $C_1(2)$, and a new partition $\Pcal_2$.  The features of
$\Pcal_2$ are that all its fibers over $C_1(2)$ have good (up to
$1 /10$) 1-block distribution, and all its fibers over $C_2(2)$
have good (up to $1 /100$) 2-block distributions.  These will not
change in the subsequent steps of the construction.

Note too that the change from $C_1(1)$, to $C_1(2)$, could have been
made arbitrarily small by choosing $N_2$ sufficiently large.

There is one problem in trying to carry out the next step
and that is, the filling in of the top relatively small portion of
the bad fibers after copying most of a good fiber.  We cannot copy
an exact good fiber because it is conceivable that no fiber with
the precise height of the bad fiber is good.  The filling in is
possible if the column heights of the previous level are relatively
prime. This was the case in step 2, because in step 1 we began with
a K-R tower heights $N_1, N_1+1$. However, Lemma \ref{4N}
does not guarantee relatively prime heights.
This is automatically the case if there is no rational
spectrum.  If there are only a finite number of rational points in
the spectrum then we could have made our original columns with
heights $LN_1,\ L(N_1+1)$ with $L$ being the highest power so that
$T^L$ is not ergodic and then worked with multiples of $L$ all the
time.  If the rational spectrum is infinite then we get an infinite
group rotation factor and this gives us the required uniform
partition without any further work.

With this point understood it is now
clear how one continues to build a sequence of
partitions $\Pcal_n$ that converge to $\Pcal$ and $C_i(k)
\rightarrow C_i$ such that the $\Pcal$-names of all fibers over
points in $C_i$ have a good (up to $1 / {10^i}$) distribution of
$i$-blocks.  This gives the uniformity of the partition
$\Pcal$ as required and establishes

\begin{prop}
Given any $\Pcal_0$ and any $\epsilon >0$ there
is a uniform partition $\Pcal$ such that
$d(\Pcal_0,\Pcal)<\epsilon$ in the $\ell_1$-metric on partitions.
\end{prop}

\br

As we have already remarked the uniform partition that we have
constructed gives us a uniquely ergodic model for the factor system
generated by this partition.  We need now a relativized version of
the construction we have just carried out.  We formulate this as
follows:

\begin{prop}
Given a uniform partition
$\Pcal$ and an arbitrary partition $\Qcal_0$ that refines
$\Pcal$, for any $\epsilon >0$ there is a uniform partition $\Qcal$
that also refines $\Pcal$ and satisfies
$$
\|\Qcal_0 - \Qcal\|_1<\epsilon.
$$
\end{prop}

Even though we write things for finite alphabets, everything makes
good sense for countable partitions as well and the arguments need
no adjusting.  However, the metric used to compare partitions
becomes important since not all metrics on $\ell_1$ are equivalent.
We use always:
$$
\|\Qcal - \overline {\Qcal}\|_1=\sum_j\ \int_X \
|1_{Q_j}-1_{\overline Q_j}| d \mu
$$
where the partitions $\Qcal$ and $\overline {\Qcal}$ are ordered
partitions into sets $\{Q_j\},\ \{\overline Q _j\}$ respectively.
We also assume that the $\sigma$-algebra generated by the partition
$\Pcal$ is nonatomic -- otherwise there is no real difference
between what we did before and what has to be done here.

We will try to follow the same proof as before.  The problem is
that when we redefine $\Qcal_0$ to $\Qcal$ we are not allowed to
change the $\Pcal$-part of the name of points.  That greatly
restricts us in the kind of names we are allowed to copy on columns
of K-R towers and it is not clear how to proceed.  The way to
overcome the difficulty is to build the K-R towers inside the
uniform algebra generated by $\Pcal$.  This being done we look,
for example, at our first tower and the first change we wish to
make in $\Qcal_0$.  We divide the fibers into a finite number of
columns according to the height and according to the $\Pcal$-name.

Next each of these is divided into subcolumns,
called $\Qcal_0$-columns, according to the
$\Qcal_0$-names of points.  If a $\Pcal$-column has some good (i.e. good
1-block distribution of $\Qcal_0$-names) $\Qcal_0$-subcolumn it can
be copied onto all the ones that are not good.  Next notice that a
$\Pcal$-column that contains not even one good $\Qcal_0$-name is
a set defined in the uniform algebra.  Therefore if these sets have
small measure then for some large enough $N$, uniformly over the
whole space, we will not encounter these bad columns too many
times.

In brief the solution is to change the nature of the uniformity. We
do not make all of the columns of the K-R tower good -- but we make
sure that the bad ones are seen infrequently, uniformly over the
whole space.  With this remark the proof of the proposition is
easily accomplished using the same nested K-R towers as before --
{\em but inside the uniform algebra}.

Finally the J-K theorem is established by constructing a refining
sequence of uniform partitions and looking at the inverse limit of
the corresponding topological spaces. Notice that if $\Qcal$
refines $\Pcal$, and both are uniform, then there is a natural
homeomorphism from $X_{\Qcal}$ onto $X_{\Pcal}$.  The way in
which the theorem is established also yields a proof of the
relative J-K theorem, Theorem \ref{rel-jk}.
\end{proof}

%

\br

Using similar methods E. Lehrer \cite{Leh} shows that
in the Jewett-Krieger theorem one can find,
for any ergodic system, a strictly
ergodic model which is topologically mixing.

\br

\section{Models for other commutative diagrams}
One can describe Theorem \ref{rel-jk} as asserting that
every diagram of ergodic systems of the form
$\Xb\to\Yb$ has a strictly ergodic model. What can we say
about more complicated commutative diagrams?
A moments reflection will show that a repeated application
of Theorem \ref{rel-jk} proves the first assertion
of the following theorem.

\begin{thm}\label{CD}
Any commutative diagram in the category of ergodic $\Z$
dynamical systems with the structure of an inverted tree,\
i.e. no portion of it looks like
\begin{equation}\label{Z>>XY}
\xymatrix
{
& \Zb \ar[dl]_\al \ar[dr]^\beta &  \\
\Xb & & \Yb
}
\end{equation}
has a strictly ergodic model.
On the other hand there exists a diagram of the form
\eqref{Z>>XY} that does not admit a strictly ergodic model.
\end{thm}

For the proof of the second assertion we need the following
theorem.

\begin{thm}\label{top-dis>meas-dis}
If $(Z,\eta,T)$ is a strictly ergodic system and
$(Z,T)\overset{\al}{\to} (X,T)$ and
$(Z,T)\overset{\beta}{\to}(Y,T)$
are topological factors such that
$\al^{-1}(U)\cap\beta^{-1}(V)\ne\emptyset$ whenever
$U\subset X$ and $V\subset Y$ are nonempty open sets,
then the measure-preserving systems $\Xb=(X,\Xcal,\mu,T)$ and
$\Yb=(Y,\Ycal,\nu,T)$ are measure-theoretically disjoint.
In particular this is the case if the systems
\XT{} and \YT{} are topologically disjoint.
\end{thm}

\begin{proof}
It suffices to show that the map $\al\times \beta:Z\to X\times Y$
is onto since this will imply that the topological system
$(X\times Y,T)$ is strictly ergodic.
We establish this by showing that the measure
$\la=(\al\times\beta)_*(\eta)$ (a joining of $\mu$ and $\nu$)
is full;\ i.e. that it assigns positive measure to
every set of the form $U\times V$ with $U$ and $V$ as in the
statement of the theorem.
In fact, since  by assumption $\eta$ is full we have
$$
\la(U\times V)=
\eta((\al\times\beta)^{-1}(U\times V))=
\eta(\al ^{-1}(U) \cap \beta^{-1}(V))>0.
$$
This completes the proof of the first assertion. The second
follows since topological disjointness of \XT{} and \YT{}
implies that  $\al\times \beta:Z\to X\times Y$ is onto.
\end{proof}

\br

\begin{proof}[Proof of Theorem Theorem \ref{CD}]
We only need to prove the last assertion.
Take $\Xb=\Yb$ to be any
nontrivial weakly mixing system, then $\Xb\times\Xb$ is
ergodic and the diagram
\begin{equation}
\xymatrix
{
& \Xb\times\Xb \ar[dl]_{p_1} \ar[dr]^{p_2} &  \\
\Xb & & \Xb
}
\end{equation}
is our counter example.
In fact if \eqref{Z>>XY} is a uniquely ergodic model in this
situation then it is easy to establish that the condition in Theorem
\ref{top-dis>meas-dis} is satisfied and we apply this theorem
to conclude that $\Xb$ is disjoint from itself.
Since in a nontrivial system
$\mu\times\mu$ and ${\gr}(\mu,{\id})$ are different
ergodic joinings, this contradiction proves our assertion.
\end{proof}

\br

\section{The Furstenberg-Weiss almost 1-1 extension theorem}

It is well known that in a topological measure space one can
have sets that are large topologically but small in the sense of the
measure. In topological dynamics when \XT is a factor of \YT
and the projection $\pi: Y \to X$ is one to one on a topologically
large set (i.e. the complement of a set of first category),
one calls \YT an {\bf almost 1-1 extension} of \XT and considers
the two systems to be very closely related. Nonetheless,
in view of the opening sentence, it is possible that the measure
theory of \YT will be quite different from the measure theory of
\XT. The following theorem realizes this possibility in an
extreme way (see \cite{FW15}).

\begin{thm}
Let \XT be a non-periodic minimal dynamical system, and let
$\pi: Y \to X$ be an extension of \XT with \YT topologically
transitive and $Y$ a compact metric space. Then there exists an
almost 1-1 minimal extension, $\ol{\pi}:(\ol{Y},T)\to (X,T)$
and a Borel subset $Y_0\subset Y$ with a Borel measurable map
$\tet: Y_0 \to \ol{Y}$ satisfying (1) $\tet T = T \tet$,
(2) $\ol{\pi}\tet =\pi$, (3) $\tet$ is 1-1 on $Y_0$, (4)
$\mu(Y_0)=1$ for any $T$-invariant measure $\mu$ on $Y$.
\end{thm}

\br

In words, one can find an almost 1-1 minimal extension of $X$
such that the measure theoretic structure is as rich as that of
an arbitrary topologically transitive extension of $X$.

An almost 1-1 extension of a minimal equicontinuous system
is called an {\bf almost automorphic} system.
The next corollary demonstrates the usefulness of this
modelling theorem. Other applications appeared e.g.
in \cite{GW35} and \cite{DL}.

\begin{cor}
Let $(X,\Xcal,\mu,T)$ be an ergodic measure preserving
transformation with infinite point spectrum defined by
$(G,\rho)$ where $G$ is a compact monothetic group
$G=\ol{\{\rho^n\}}_{n\in \Z}$. Then there is an almost
1-1 minimal extension of $(G,\rho)$ (i.e. a minimal almost
automorphic system), $(\tilde Z,\sig)$ and an invariant
measure $\nu$ on $Z$ such that $(Z,\sig,\nu)$ is isomorphic
to $(X,\Xcal,\mu,T)$.
\end{cor}

\br

\section{Cantor minimal representations}

A {\em Cantor minimal dynamical system} is
a minimal topological system \XT where $X$ is the Cantor set.
Two Cantor minimal systems \XT and $(Y,S)$ are called
{\em orbit equivalent\/} (OE)
if there exists a homeomorphism $F: X\to Y$
such that $F(\OT(x))=\OS(Fx)$ for every $x\in X$.
Equivalently: there are functions
$n:X\to \Z$ and $m: X \to \Z$ such that for every $x\in X$
$F(Tx)=S^{n(x)}(Fx)$ and $F(T^{m(x)})=S(Fx)$.
An old result of M. Boyle implies that the requirement that,
say, the function $n(x)$ be continuous already implies that
the two systems are {\em flip conjugate\/}; i.e. $(Y,S)$
is isomorphic either to \XT or to $(X,T^{-1})$. However,
if we require that both $n(x)$ and $m(x)$ have at most
one point of discontinuity we get the new and, as it turns out,
useful notion of {\em strong orbit equivalence\/} (SOE).
A complete characterization of both OE and SOE of Cantor minimal systems
was obtained by Giordano Putnam and Skau \cite{GPS}
in terms of an algebraic invariant of Cantor minimal systems called
the {\em dimension group}. (See \cite{Glasner} for
a review of these results.)

We conclude this section with
the following remarkable theorems, due to N. Ormes \cite{Ormes},
which simultaneously generalize the theorems of Jewett and Krieger
and a theorem of Downarowicz \cite{Dow} which,
given any Choquet simplex $Q$, provides
a Cantor minimal system \XT with $M_T(X)$ affinely homeomorphic with $Q$.
(See also Downarowitcz and Serafin \cite{DS},
and Boyle and Downarowicz \cite{BD}.)

\begin{thm}\label{ort}
\begin{enumerate}
\item
Let $(\Om,\Bcal,\nu,S)$ be an ergodic, non-atomic, probability measure
preserving, dynamical system. Let $(X,T)$ be a Cantor minimal
system such that whenever $\exp(2\pi i/p)$
is a (topological) eigenvalue of
\XT for some $p\in \N$ it is also a (measurable) eigenvalue
of $(\Om,\Bcal,\nu,S)$.
Let $\mu$ be any element of the set of extreme points
of $M_T(X)$. Then, there exists
a homeomorphism $T':X \to X$ such that (i) $T$ and $T'$ are
strong orbit equivalent, (ii) $(\Om,\Bcal,\nu,S)$ and $(X,\Xcal,\mu,T')$
are isomorphic as measure preserving dynamical systems.
\item
Let $(\Om,\Bcal,\nu,S)$ be an ergodic, non-atomic, probability measure
preserving, dynamical system. Let $(X,T)$ be a Cantor minimal
system and $\mu$ any element of the set of extreme points
of $M_T(X)$. Then, there exists
a homeomorphism $T':X \to X$ such that (i) $T$ and $T'$ are
orbit equivalent, (ii) $(\Om,\Bcal,\nu,S)$ and $(X,\Xcal,\mu,T')$
are isomorphic as measure preserving dynamical systems.
\item
Let $(\Om,\Bcal,\nu,S)$ be an ergodic, non-atomic, probability measure
preserving dynamical system. Let $Q$ be any Choquet simplex and
$q$ an extreme point of $Q$. Then there exists a
Cantor minimal system \XT and an affine homeomorphism $\phi: Q \to M_T(X)$
such that, with $\mu=\phi(q)$,
$(\Om,\Bcal,\nu,S)$ and $(X,\Xcal,\mu,T)$
are isomorphic as measure preserving dynamical systems.
\end{enumerate}
\end{thm}

\br

\section{Other related theorems}
Let us mention a few more striking representation results.

For the first one recall that a topological dynamical system
$(X,T)$ is said to be \textbf{prime} if it has no non-trivial
factors. A similar definition can be given for measure preserving
systems. There it is easy to see that a prime system $(X,\Xcal,\mu,T)$
must have zero entropy. It follows  from a construction in \cite{SW}
that the same holds for topological entropy, namely any system $(X,T)$
with positive topological entropy has non-trivial factors. In
\cite{Wei36} it is shown that any ergodic zero entropy
dynamical system has a minimal model $(X,T)$ with the property
that any pair of points $(u,v)$ not on the same orbit has a dense
orbit in $X \times X$. Such minimal systems are necessarily prime,
and thus we have the following result:

\begin{thm}\label{prime}
An ergodic dynamical system has a topological, minimal,
prime model iff it has zero entropy.
\end{thm}

The second theorem,
due to Glasner and Weiss \cite{GW35},
treats the positive entropy systems.

\begin{thm}
An ergodic dynamical system has a strictly ergodic,
UPE model iff it has positive entropy.
\end{thm}

\br

We also have the following surprising result which is due
to Weiss \cite{Wei35}.

\begin{thm}
There exists a minimal metric dynamical system $(X,T)$
with the property that for every ergodic probability
measure preserving system $(\Om,\Bcal,\mu,S)$ there
exists a $T$-invariant Borel probability measure
$\nu$ on $X$ such that the systems
$(\Om,\Bcal,\mu,S)$ and $(X,\Xcal,\nu,T)$ are isomorphic.
\end{thm}

\br

In \cite{Lis2} E. Lindenstrauss proves the following:
\begin{thm}\label{distal-model}
Every ergodic measure distal $\Z$-system
$\Xb=(X,\Xcal,\mu,T)$ can be
represented as a minimal topologically distal
system $(X,T,\mu)$ with $\mu\in M_T^{\erg}(X)$.
\end{thm}
This topological model need not, in general,
be uniquely ergodic. In other words there are measure
distal systems for which no uniquely ergodic topologically distal
model exists.

\begin{prop}
\begin{enumerate}
\item
There exists an ergodic non-Kronecker measure distal system
$(\Omega,\Fcal,m,T)$ with nontrivial maximal
Kronecker factor
$(\Om_0,\Fcal_0,m_0,T)$ such that (i) the extension
$(\Omega,\Fcal,m,T)\to (\Om_0,\Fcal_0,m_0,T)$ is finite to one a.e.
and (ii) every nontrivial factor map of $(\Om_0,\Fcal_0,m_0,T)$
is finite to one.
\item
A system $(\Omega,\Fcal,m,T)$ as in part 1 does not admit
a topologically distal strictly ergodic model.
\end{enumerate}
\end{prop}

\begin{proof}
1.\ Irrational rotations of the circle as well as
adding machines are examples of Kronecker systems satisfying
condition (ii).
There are several constructions in the literature of ergodic,
non-Kronecker, measure distal, two point extensions of these
Kronecker systems.
A well known explicit example is the strictly ergodic Morse
minimal system.

\br

2.\
Assume to the contrary that $(X,\mu,T)$ is a distal
strictly ergodic model for $(\Omega,\Fcal,m,T)$.
Let $(Z,T)$ be the maximal equicontinuous factor of
$(X,T)$ and let $\eta$ be the unique invariant probability
measure on $Z$. Since by assumption $(X,\mu,T)$ is
not Kronecker it follows that $\pi: X \to Z$ is not
one to one. By Furstenberg's structure theorem for minimal distal
systems $(Z,T)$ is nontrivial and moreover there exists an
intermediate extension
$X \to Y \overset{\sigma}{\to}Z$ such that $\sigma$ is
an isometric extension. A well known construction implies
the existence of a minimal group extension
$\rho:(\tilde Y,T) \to (Z,T)$, with compact fiber group $K$,
such that the following diagram is commutative
(see Section \ref{Sec-distal} above).
We denote by $\nu$ the unique invariant measure on $Y$
(the image of $\mu$) and let $\tilde \nu$ be an ergodic
measure on $\tilde Y$ which projects onto $\nu$.
The dotted arrows denote measure theoretic factor maps.

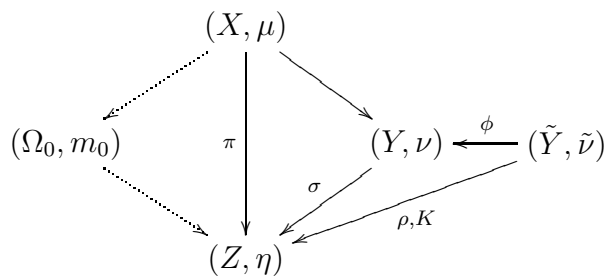
\begin{equation*}
\xymatrix
{
& (X,\mu)\ar@{.>}[dl]\ar[dd]_{\pi} \ar[dr]&  & &\\
(\Om_0,m_0) \ar@{.>}[dr] & &(Y,\nu)\ar[dl]_{\sigma} &
(\tilde Y,\tilde\nu)\ar[l]_{\phi}
\ar[dll]^{\rho, K}\\
& (Z,\eta)& & &
}
\end{equation*}

\br

Next form the measure
$ \theta = \int_K R_k\tilde\nu\, dm_K,$
where $m_K$ is Haar measure on $K$ and for each $k\in K$,
$R_k$ denotes right translation by $k$ on $\tilde Y$
(an automorphism of the system $(\tilde Y,T)$).
We still have $\phi(\theta)=\nu$.

A well known theorem in topological dynamics (see \cite{SS}) implies
that a minimal distal finite to one extension of a minimal equicontinuous
system is again equicontinuous and since $(Z,T)$ is the maximal
equicontinuous factor of $(X,T)$ we conclude that
the extension $\sigma: Y \to Z$ is not finite to one.
Now the fibers of the extension $\sigma$ are homeomorphic
to a homogeneous space $K/H$, where $H$ is a closed subgroup of $K$.
Considering the measure disintegration $\theta = \int_Z \theta_z\,
d\eta(z)$ of $\theta$ over $\eta$ and its projection
$\nu = \int_Z \nu_z\, d\eta(z)$, the disintegration of $\nu$ over
$\eta$, we see that a.e. $\theta_z \equiv m_K$ and
$\nu_z \equiv m_{K/H}$. Since $K/H$ is infinite we conclude
that the {\em measure theoretical extension}
$\sigma: (Y,\nu) \to (Z,\eta)$
is not finite to one. However considering the dotted part
of the diagram we arrive at the opposite conclusion.
This conflict concludes the proof of the proposition.
\end{proof}

\br

In \cite{OW} Ornstein and Weiss introduced the notion
of tightness for measure preserving systems and the analogous
notion of mean distality for topological systems.

\begin{defn}
Let \XT be a topological system.
\begin{enumerate}
  \item
  A pair $(x,y)$ in $X\times X$ is {\bf mean proximal\/}
  if for some (hence any) compatible metric $d$
  $$
  \limsup_{n\to\infty}\frac{1}{2n+1}\sum_{i=-n}^{n}
  d(T^i x, T^i y) = 0.
  $$
  If this $\limsup$ is positive the pair is called
  {\bf mean distal\/}.
  \item
  The system \XT is {\bf mean distal\/} if every pair with $x\ne y$
  is mean distal.
  \item
  Given a $T$-invariant probability measure $\mu$ on $X$,
  the triple $(X,\mu,T)$ is called {\bf tight\/} if there is
  a $\mu$-conull set $X_0\subset X$ such that every pair of
  distinct points $(x,y)$ in $X_0\times X_0$ is mean distal.
\end{enumerate}
\end{defn}

Ornstein and Weiss show that tightness is in fact a property
of the measure preserving system $(X,\mu,T)$ (i.e. if the
measure system $(X,\Xcal,\mu,T)$ admits one tight model then
every topological model is tight). They obtain the following results.

\begin{thm}
\mbox{}
\begin{enumerate}
  \item
  If the entropy of $(X,\mu,T)$ is positive and finite
  then $(X,\mu,T)$ is not tight.
  \item
  There exist strictly ergodic non-tight systems
  with zero entropy.
\end{enumerate}
\end{thm}

Surprisingly the proof in \cite{OW} of the
non-tightness of a positive entropy system does not
work in the case when the entropy is infinite which is
still open.

J. King gave an example of a tight system with a non-tight
factor. Following this he and Weiss \cite{OW} established
the following result. Note that this theorem
implies that tightness and mean distality are not
preserved by factors.

\begin{thm}
  If  $(X,\Xcal,\mu,T)$ is ergodic with zero entropy
  then there exists a mean-distal system $(Y,\nu,S)$ which
  admits $(X,\Xcal,\mu,T)$ as a factor.
\end{thm}

\br


\end{document}